\documentclass{amsart}
\usepackage{amsfonts,amssymb,amsmath,amsthm}
\usepackage{url}
\usepackage{enumerate}

\usepackage{amsmath,amsfonts,amsthm,url,color,amssymb}
\usepackage{graphicx}
\usepackage{algpseudocode, algorithm}
\newtheorem{theorem}{Theorem}[section]
\newtheorem{remark}{Remark}

\newtheorem{proposition}[theorem]{Proposition}
\newtheorem{lemma}[theorem]{Lemma}
\newtheorem{define}[theorem]{Definition}
\newtheorem{example}{Example}
\newtheorem{corollary}[theorem]{Corollary}
\newcommand{\rmv}[1]{}
\numberwithin{equation}{section}
\newcommand{\I}{\mathcal I}
\newcommand{\ou}{\mathcal O}
\newcommand{\cc}{\mathrm{Cyc}}
\newcommand{\sym}{\mathrm{Sym}}
\newcommand{\C}{\mathcal C}
\newcommand{\mq}{m_{\alpha, q}(x)}

\newcommand{\G}{\mathbb G}

\newcommand{\GL}{\mathrm{GL}}
\newcommand{\PGL}{\mathrm{PGL}}

\newcommand{\F}{\mathbb{F}}
\newcommand{\fq}{\mathbb{F}_q}
\newcommand{\ord}{\mathrm{ord}}
\newcommand{\oord}{\mathrm{Ord}}

\author[L. Reis]{Lucas Reis}
\address{
Departamento de Matem\'{a}tica\\
Universidade Federal de Minas Gerais\\
Belo Horizonte MG, Brazil}
\email{lucasreismat@gmail.com}
\thanks{This work was conducted during a visit of the first author to Carleton University, supported by CAPES-PDSE (process - 88881.134747/2016-01). The research of Qiang Wang is partially supported by NSERC of Canada.}

\author[Q. Wang]{Qiang Wang}
\address{
School of Mathematics and Statistics\\
Carleton University\\
Ottawa ON, Canada}
\email{wang@math.carleton.ca}

\begin{document}

\title[The dynamics of permutations on irreducible polynomials]{The dynamics of permutations on irreducible polynomials}

\begin{abstract}
We study degree preserving maps over the set of irreducible polynomials over a finite field. In particular, we show that every permutation of the set of irreducible polynomials of degree $k$ over $\mathbb{F}_q$ is induced by an action from a permutation polynomial of $\mathbb{F}_{q^k}$ with coefficients in $\mathbb{F}_q$. The dynamics of these permutations of irreducible polynomials of degree $k$ over $\mathbb{F}_q$, such as fixed points and cycle lengths,  are studied. As an application, we also generate irreducible polynomials of the same degree by an iterative method.  
\end{abstract}
\maketitle

\section{Introduction}
Let $\F_q$ be the finite field with $q$ elements, where $q$ is a power of a prime $p$ and let $\overline{\F}_q$ be the algebraic closure of $\F_q$.  For any $\alpha\in  \overline{\F}_q$, we let $m_{\alpha}(x)$ denote the minimal polynomial of $\alpha$ over $\F_q$. A polynomial $F(x) \in \F_q[x]$ naturally induces a function in the set $\overline{\F}_q$: namely, we have the evaluation map $F:\overline{\F}_q\to \overline{\F}_q$ with $c\mapsto F(c)$. For each positive integer $k$, let $\I_k$ be the set of monic irreducible polynomials over $\F_q$ of degree $k$ and $\I=\cup_{k\ge 1}\I_k$ be the set of monic irreducible polynomials over $\F_q$.  It turns out that any polynomial $F(x) \in\F_q[x]$ induces a map on the set $\I$; namely,  for any $f\in \I$, let $\alpha\in \overline{\F}_q$ be any root of $f$, we define $F\diamond f$ as the minimal polynomial of $\beta=F(\alpha)$, i.e., 
$$F\diamond f:= m_{F(\alpha)}(x).$$

Of course, $F\diamond f$ is again in $\I$. Since $F(x)\in \F_q[x]$, this composition operation $\diamond$ is well defined: if $\alpha_0$ is another root of $f$, then $\alpha_0=\alpha^{q^i}$ is a conjugate of $\alpha$ over $\F_q$ and so $F(\alpha_0)=F(\alpha^{q^i})=F(\alpha)^{q^i}$ is a conjugate of $F(\alpha)$, i.e., $m_{F(\alpha)}=m_{F(\alpha_0)}$. In particular, the polynomial $F$ yields a dynamics on the set $\I$: for $f\in \I$, one may consider the orbit $\{f, F\diamond f,  F\diamond (F\diamond f), \ldots\}$ of $f$ by $F$. The study of this kind of dynamics was initiated by Vivaldi~\cite{V92} and Batra and Morton~\cite{BM94-1, BM94-2}. In \cite{BM94-1, BM94-2}, the authors explored this dynamics over $\I$ for some special classes of linearized polynomials and, in particular, they show the existence of infinitely many fixed points, i.e., orbits of length one. This study was extended in \cite{M97}, where Morton showed that, for $F(x) =x^q-ax, a\in \F_q^*$, the dynamics induced by $F$ on $\I$ yields infinitely many periodic points with period $\pi_F(f)=n$, for any positive integer $n$. In \cite{CH00}, Cohen and Hachenberger  extended this last result to any linearized polynomial $F(x) =\sum_{i=0}^{s}a_ix^{q^i}$ not of the form $ax^{q^s}$ for any $a\in \F_q^*$. 

In general, if $f\in \I_k$ and $\alpha\in \F_{q^k}$ is any root of $f$, then $F(\alpha)$ is in $\F_{q^k}$ and  its minimal polynomial is of degree $d$, where $d$ is a divisor of $k$. In particular, we have $\deg(F\diamond f)\le \deg(f)$ for any $F\in \F_q[x]$ and $f\in \I$.  It is then natural to study the extremal case, i.e., $\deg(F\diamond f)=\deg(f)$ for any $f\in \I$. In other words, we are interested in any polynomial $F$ for which $f\mapsto F\diamond f$ is a \emph{degree preserving} map. We shall prove in Theorem~\ref{thm:degree-preserve}  that the class of polynomials $F$ inducing a degree preserving map on $\I$ is extremely small, namely, the class of polynomials $F(x)=ax^{p^h}+b$ with $a, b\in \F_q$, $a \neq 0$, and $h\ge 0$.


Nevertheless, we may obtain many examples of polynomials $F$ for which $f\mapsto F\diamond f$ is \emph{locally} degree preserving. We recall that a polynomial $f\in \F_{q^k}[x]$ is called a permutation polynomial of $\F_{q^k}$ if the evaluation mapping is a permutation of $\F_{q^k}$. We prove that, if $F(x) \in \F_q[x]$ is a permutation polynomial of $\F_{q^k}$, then   $F\diamond f \in \I_k$ for any $f \in \I_k$. This means that $\deg(F\diamond f) = \deg(f)$ for all $f \in \I_k$. Moreover,  the map $f\mapsto F\diamond f$ from $\I_k$ to the set $\I_k$ is a permutation; for more details, see Theorem~\ref{thm:local}.

Our paper mostly concentrates on permutations of the set $\I_k$ of monic irreducible polynomials of degree $k$. 
 We prove in Theorem~\ref{thm:onto} that every permutation $\sigma$ of the set $\I_k$ can be represented as the map induced by a permutation polynomial $P(x)$ of $\F_{q^k}$ with all coefficients in $\F_q$, i.e,, $\sigma(f) = P\diamond f$ for any $f\in \I_k$. To study the dynamics of $P(x)$ on $\I_k$, we find that it is more convenient to study an associated operation $P\ast f = Q\diamond f$, where $Q$ is the compositional inverse of $P$. In particular, we show that

\begin{equation*}P\ast f=\gcd(f(P(x)), x^{q^k}-x).\end{equation*}

In Section 5, we provide some general results on the fix points of the compositions $P\ast f$. We prove that
$P\ast f = f$ if and only if $f$ is a divisor of $\prod_{i=0}^{k-1} (x^{q^i} -P(x))$. Through this characterization, we derive the number of fixed points under this composition operation. In Section 6, we provide connections between the dynamics of the 
evaluation map  $c\mapsto P(c)$ over $\F_{q^k}$ (restricted on the set $\C_k$ of elements of degree $k$ over $\F_q$)  
 and the map $f \mapsto P\ast f$ over $\I_k$.
  Finally, we propose to use this operation to iterately generate irreducible polynomials of the same degree using any given irreducible polynomial. In some cases, we show that we can generate at least a proportion $1/k$ of all the irreducible polynomials of degree $k$, using any one irreducible polynomial of the degree $k$. \\


\section{Preliminaries}
In this section, we provide a background material that is used along the paper. Most of the following results are standard in the theory of finite fields and, for this reason, we skip some details.

\subsection{Multiplicative order in finite fields}\label{subsec:mult}

 For $\alpha\in  \overline{\F}_q^{*}$, the multiplicative order of $\alpha$ is defined by $\ord(\alpha):=\min\{d>0\,|\, \alpha^d=1\}$. The degree $\deg(\alpha)$ of $\alpha\in\overline{\F}_q$ over $\F_q$ is defined as the degree of the minimial polynomial of $\alpha$ over $\F_q$.  
If $a, n$ are positive integers such that $\gcd(a, n)=1$, then $\ord_{n}a:=\min\{r>0\,|\, a^r\equiv 1\pmod n\}$. In the following theorem, we present some well-known results on the multiplicative order of elements in finite fields.

\begin{lemma}\label{thm:mult}
Let $\alpha\in \overline{\F}_q^*$ be an element of multiplicative order $\ord(\alpha)=e$. The following hold:
\begin{enumerate}[(i)]
\item $\deg(\alpha)=\ord_eq$;
\item if $\beta=\alpha^s$, then $\ord(\beta)=\frac{e}{\gcd(e, s)}$.
\end{enumerate}
For any positive integer $E$ that is relatively prime with $q$, there exist $\varphi(E)$ elements $\alpha\in \overline{\F}_q^*$ such that $\ord(\alpha)=E$, where $\varphi$ is the \emph{Euler Function}. In addition, such elements are all in $\F_{q^k}$, where $k=\ord_Eq$.
\end{lemma}

If $f\in \I_k$ with $f(x)\ne x$, then the order $\ord(f)$ of $f$ is the order of any (hence all) of its roots.  An element $\alpha\in \F_{q^k}^{*}$ is \emph{primitive} if $\ord(\alpha)=q^k-1$. Additionally, $f\in \I_k$ is a \emph{primitive polynomial} if any (hence all) of its roots is a primitive element, i.e., $\ord(f)=q^k-1$. In particular, Lemma~\ref{thm:mult} shows that there exist $\varphi(q^k-1)>0$ primitive elements in $\F_{q^k}$.

\subsection{Linearized polynomials and the $\F_q$-order}\label{subsec:add}
A polynomial of the form $L(x)=\sum_{i=0}^{m}a_ix^{q^i}\in \overline{\F}_q[x]$ is called \emph{linearized}. Given $f\in \overline{\F}_q[x]$ with $f(x)=\sum_{i=0}^{n}a_ix^i$, we can associate to $f$ the polynomial $L_f(x)=\sum_{i=0}^{n}a_ix^{q^i}$. In this case, $L_f$ is the \emph{$q$-associate} of $f$.  Of course, any linearized polynomial is the $q$-associate of some polynomial in $\overline{\F}_q[x]$. In the following lemma, we show that linearized polynomials with coefficients in $\F_q$ behave well through basic operations.

\begin{lemma}\label{lem:linearized-properties}
Let $f, g\in \F_{q}[x]$. The following hold:
\begin{enumerate}[(i)]
\item $L_f(x)+L_g(x)=L_{f+g}(x)$,
\item $L_g(L_f(x)=L_{fg}(x)=L_{gf}(x)$,
\item $\gcd(L_f, L_g)=L_{\gcd(f, g)}$. 
\end{enumerate}
\end{lemma}
\begin{proof}
Items (i) and (ii) follow by direct calculations and item (iii) is proved in Section 3.4 of \cite{LN}. 
\end{proof}

For a polynomial $f\in \F_q[x]$ with $f(x)=\sum_{i=0}^{m}a_ix^i$ and an element $\alpha\in \overline{\F}_q$, we have  $L_f(\alpha)=\sum_{i=0}^ma_i\alpha^{q^i}$. 
We observe that, if $\alpha\in \F_{q^k}$, $L_{x^{k}-1}(\alpha)=\alpha^{q^k}-\alpha=0$. We define $$\mathcal A_{\alpha}=\{g \in \F_q[x]\, |\, L_g(\alpha)=0\}.$$ From the previous lemma, $\mathcal A_{\alpha}$ is an ideal of $\F_q[x]$ and we observe that $x^k-1$ is in this ideal. In particular, $\mathcal A_{\alpha}$ is generated by a non-zero polynomial $\mq$, which we can suppose to be monic. 

\begin{define}The polynomial $\mq$ is defined as the \emph{$\F_q$-order} of $\alpha$. Namely, $\mq$ is a polynomial in $\F_q[x]$ with the lowest degree such that $\alpha$ is a root of its $q$-associate. \end{define}

For instance the element $0$ has $\F_q$-order $m_{0, q}(x)=1$ and, for any $c\in \F_{q}^*$, $m_{c, q}(x)=x-1$. For any $\alpha\in \overline{\F}_q$ and any $f\in \F_q[x]$,
by Lemma~\ref{lem:linearized-properties}, we have that $L_f(\alpha)=0$ if and only if $f$ is divisible by $\mq$. The polynomial $\mq$ works as an ``additive'' order, in duality to the multiplicative order over finite fields. It is clear that $\alpha$ is in $\F_{q^k}$ if and only if $0=\alpha^{q^k}-\alpha=L_{x^k-1}(\alpha)$, i.e., $\mq$ divides $x^k-1$. Also, $\alpha$ and any of its conjugates $\alpha^{q^i}$ have the same $\F_q$-order. 
This motivates us to introduce the following definition.
\begin{define}
For $f\in \F_q[x]$ an irreducible polynomial and let $\alpha$ be one of its roots, the $\F_q$-order of $f$ is defined as the $\F_q$-order of $\alpha$. We write $\oord(f)=\mq$.  
\end{define}


The elements $\alpha\in \F_{q^k}$ for which $\mq=x^k-1$ are called {\it normal} over $\F_q$. Normal elements work as ``additive generators'' (in the linearized sense) of the additive group $\F_{q^k}$, in duality to the primitive elements for the multiplicative group $\F_{q^k}^*$. Following the definition of primitive polynomials, an element $f\in \I_k$ is a \emph{normal polynomial} if $\oord(f)=x^k-1$. 

In the additive-multiplicative order analogies, Lemma~\ref{thm:mult} can be translated to $\F_q$-order with a suitable change of functions. 

\begin{define}
Let $f, g\in \F_q[x]$.
\begin{enumerate}[(i)]
\item The \emph{Euler Phi function} for polynomials over $\F_q$ is $$\Phi_q(f)=\left |\left(\frac{\F_q[x]}{\langle f\rangle}\right)^{*}\right |,$$ where $\langle f\rangle$ is the ideal generated by $f$ in $\F_q[x]$;
\item If $\gcd(f, g)=1$, then $\ou(f, g):=\min\{k>0\,|\, f^k\equiv 1\pmod g\}$ is the order of $f$ modulo $g$.
\end{enumerate}
\end{define}

The function $\Phi_q$ is multiplicative. Also, $\Phi_q(g^s)=q^{(s-1)d}(q^d-1)$ if $g$ is an irreducible polynomial of degree $d$ and $s$ is a positive integer: one may compare with $\varphi(r^s)=r^{s-1}(r-1)$ if $r$ is a prime number. It is straightforward to check that $\ou(f, g)$ divides $\Phi_q(g)$.

\begin{theorem}\label{thm:add}
Let $\alpha\in \overline{\F}_q$ and the $\F_q$-order of $\alpha$ be $\mq=h(x)$. The following hold:
\begin{enumerate}[(i)]
\item $\deg(\alpha)=\ou(x, h(x))$;
\item if $\beta=L_g(\alpha)$, then $\beta$ has $\F_q$-order $m_{\beta, q}(x)=\frac{h(x)}{\gcd(h(x), g(x))}$.
\end{enumerate}
In addition, for any polynomial $H$ relatively prime with $x$, there exist $\Phi_q(H)$ elements $\alpha\in \overline{\F}_q$ such that $\mq=H$.
\end{theorem}

\begin{proof}
\begin{enumerate}[(i)]
\item Observe that $\deg(\alpha)$ is the least positive integer $k$ such that $\alpha\in \F_{q^k}$. Also, for any positive integer $d$, we have that $\alpha\in \F_{q^d}$ if and only if $$L_{x^d-1}(\alpha)=\alpha^{q^d}-\alpha=0,$$ that is, $\mq=h(x)$ divides $x^d-1$. Now, the result follows from definition of $\ou(x, h)$.

\item This item follows by direct calculations.
\end{enumerate}
For the proof of the last statement, see Theorem 11 of~\cite{O34}.
\end{proof}

\subsection{Permutation polynomials}\label{subsec:permut}
Here we present some well-known classes of permutations of finite fields.

\begin{enumerate}[$\bullet$]
\item \emph{Monomials}: The polynomial $x^n$ is a permutation of $\F_{q^k}$ if and only if $$\gcd(n, q^k-1)=1.$$

\item \emph{Linearized}: If $g(x) \in \F_{q}[x]$ with $g(x) =\sum_{i=0}^{k-1}a_ix^i$ and $L_g(x)=\sum_{i=0}^{k-1}a_ix^{q^i}$ is the \emph{$q$-associate of $f$}, then $L_g$ is a permutation of $\F_{q^k}$ if and only if $$\gcd(g(x), x^k-1)=1.$$

\item \emph{M\"{o}bius}: Let $\GL_2(\F_q)$ be the group of all invertible $2\times 2$ matrices over $\F_q$. Given $[A]\in \GL_2(\F_q)$ with $A=\left(\begin{matrix}
a&b\\
c&d
\end{matrix}\right)
$, we set $$\tau_A:\F_{q^k}\to \F_{q^k},$$ where $\tau_A(z)=\frac{az+b}{cz+d}$ if $c=0$ or $c\ne 0$ and $z\ne -d/c$, and $\tau_A(-d/c)=a/c$ if $c\ne 0$. Then $\tau_A$ is a permutation of $\F_{q^k}$ for any $k\ge 1$.
The map $\tau_A(z)$ admits a polynomial representation; for any $z\in \F_{q^k}$, $\tau_A(z)=(az+b)d^{-1}$ if $c=0$ and $\tau_A(z)=(az+b)\left[(cz+d)^{q^k-2}+\varepsilon\cdot \frac{z^{q^k}-z}{z+d/c}\right]$ if $c\ne 0$, where $\varepsilon=\frac{a}{\det A}$. 
\end{enumerate}

\subsection{Notations} We want to emphasize that the following notations and consequences are frequently employed in this paper. 

\begin{enumerate}[$\bullet$]

\item If $a, n$ are positive integers such that $\gcd(a, n)=1$, then $\ord_{n}a:=\min\{r>0\,|\, a^r\equiv 1\pmod n\}$.

\item $\overline{\F}_q$ denotes the algebraic closure of $\F_q$.

\item For $\alpha\in  \overline{\F}_q^{*}$, $\ord(\alpha):=\min\{d>0\,|\, \alpha^d=1\}$ is the multiplicative order of $\alpha$.

\item $\I_k$ denotes  the set of irreducible monic polynomials of degree $k$ over $\F_q$.


\item For $\alpha\in  \overline{\F}_q$, $m_{\alpha}(x)$ is the minimal polynomial of $\alpha$ over $\F_q$. 

\item $\deg(\alpha):= \deg(m_{\alpha})=\min\{s>0\,|\, \alpha\in \F_{q^s}\}$.

\item  For $\alpha\in  \overline{\F}_q$, the $\fq$-order $\mq$ is the polynomial of the lowest degree polynomial such that $\alpha$ is a root of its $q$-associate.

\item   For $\alpha\in  \overline{\F}_q$,  $\alpha$ is in $\F_{q^k}$ if and only if  $\mq$ divides $x^k-1$.

\item $\C_k$ denotes the set of elements $\alpha \in \overline{\F}_q$ such that $\deg(\alpha)=k$, i.e., $m_{\alpha}(x)\in \I_k$.
\item $\G_k:=\{P\in \F_q[x]\,|\, P\; \text{is a permutation polynomial of}\,\, \F_{q^k}\, \text{and}\, \deg(P)<q^k\}.$

\end{enumerate}

\section{On degree preserving maps}
In this section, we provide the proofs of Theorems~\ref{thm:degree-preserve} and~\ref{thm:local}.  We start with recalling some notations: for $\alpha\in \overline{\F}_q$,  $\deg(\alpha):= \deg(m_{\alpha})=\min\{s>0\,|\, \alpha\in \F_{q^s}\}$. Moreover, for each positive integer $k$,  $\C_k$ is the set of elements $\alpha \in \overline{\F}_q$ such that $\deg(\alpha)=k$, i.e., $m_{\alpha}\in \I_k$. In order to prove Theorem~\ref{thm:degree-preserve}, the following proposition is crucial.

\begin{proposition}\label{prop:degree-preserving}
Suppose that $F\in \F_q[x]$ is a polynomial of degree $d\ge 1$ such that its induced map $f\mapsto F\diamond f$ on $\I$ preserves the degree of the elements in $\I$. Then, for each $\alpha\in \overline{\F}_q$, the equation $F(x)=\alpha$ has exactly one solution $\gamma\in \overline{\F}_q$ with multiplicity $d$ and, in fact, $\deg(\gamma)=\deg(\alpha)$. In addition, for any positive integer $k$, the evaluation map induced by $F$ on $\F_{q^k}$ is a permutation and, in particular, the evaluation map induced by $F$ over the field $\overline{\F}_q$ is a permutation.
\end{proposition}
\begin{proof}
Let $j$ be any positive integer and, for each $\alpha\in \C_j$, let $C(\alpha, F)\subseteq \overline{\F}_q$ be the set of (distinct) solutions of $F(x)=\alpha$, $|C(\alpha, F)|\ge 1$. For each $\gamma\in C(\alpha, F)$, we have $F(\gamma)=\alpha$ hence $F\diamond m_{\gamma}(x)=m_{\alpha}(x)$. Since $F$ preserves degree, it follows that $\deg(\gamma)=\deg(\alpha)=j$, hence $\gamma\in \C_j$. In particular, for each $\alpha\in \C_j$, $C(\alpha, F)\subseteq \C_j$. It is straightforward to check that the sets $C(\alpha, F)$ are pairwise disjoint for distinct $\alpha$'s. Since $\bigcup_{\alpha\in \C_j}C(\alpha, F)\subseteq \C_j$ is a disjoint union of nonempty sets, it follows that
$$|\C_j|\le \sum_{\alpha\in \C_j}|C(\alpha, F)|\le |\C_j|,$$
hence $|C(\alpha, F)|=1$ for any $\alpha$. In other words, $F(x)-\alpha$ has one root $\gamma\in \C_j$ with multiplicity $d$. In particular, since $\F_{q^k}=\bigcup_{j|k}\C_j$ one can see that, for each $\gamma\in \F_{q^k}$, the equation $F(x)=\gamma$ has exactly one solution over $\overline{\F}_q$ and this solution lies in $\F_{q^k}$. Therefore, the evaluation map $a\mapsto F(a)$ on $\F_{q^k}$ (a finite set) is onto and so is a permutation.
\end{proof}

As follows, we provide a complete characterization of the polynomials $F$ satisfying the properties given in Proposition~\ref{prop:degree-preserving}.

\begin{lemma}\label{lem:trivials}
Suppose that $F\in \F_q[x]$ is a polynomial of degree $d\ge 1$ such that, for each $\alpha\in \overline{\F}_q$, the $F(x)-\alpha$ has a unique root in $\overline{\F}_q$ with multiplicity $d$. Then $F(x)=ax^{p^h}+b$ for some $a, b\in \F_q$ and some $h\ge 0$.
\end{lemma}
\begin{proof}
Write $F(x)=\sum_{i=0}^{d}a_ix^i$ and let $p^h$ be the greatest power of $p$ that divides each index $i$ for which $a_i\ne 0$. Hence $F=G^{p^h}$, where $G=\sum_{i=0}^sa_i^\prime x^i\in \F_q[x]$ is such that there exists at least one positive integer $1\le i\le s$ for which $a_i^\prime \ne 0$ and $i$ is not divisible by $p$. In particular $G'(x)$, the formal derivative of $G$, is not the zero polynomial. One can easily see that $G$ also satisfies the required properties of our statement. We shall prove that $G$ has degree $s=1$. For this, suppose that $s>1$ and, for each $\alpha\in \overline{\F}_q$, let $\gamma(\alpha)$ be the only root of $G(x)=\alpha$, hence $\gamma(\alpha)$ has multiplicity $s\ge 2$. In particular, $\gamma(\alpha)$ is a root of $\gcd(G(x)-\alpha, G'(x))$. This shows that $G'(x)$ vanishes at each element $\gamma(\alpha)$ with $\alpha\in \overline{\F}_q$. Of course, the set of elements $\gamma(\alpha)$ is infinite and so $G'$ is the zero polynomial, a contradiction. Therefore, $G$ has degree $s=1$ and so $G=Ax+B$, hence $F(x) =ax^{p^h}+b$ with $a=A^{p^h}$ and $b=B^{p^h}$.
\end{proof}

We observe that, if $F(x)=ax^{p^h}+b$ is a polynomial of degree $p^h\ge 1$ (hence $a\ne 0$), then $\deg(F\diamond f)=\deg(f)$ for any monic irreducible $f$. For this, suppose that $\deg(f)=k$ and let $\alpha\in \F_{q^k}$ be any root of $f$. It is straightforward to check that $F$ permutes the whole field $\overline{\F}_q$. In particular, since the compositions $F^{(n)}(\alpha)$ are in $\F_{q^k}$ (which is a finite set), there exists a positive integer $m$ such that $F^{(m)}(\alpha)=\alpha$. Therefore, if we set $f_0=f$ and, for $1\le i\le m$, $f_i=F\diamond f_{i-1}$, one has $f_m=F^{(m)} \diamond f=m_{F^{(m)}(\alpha)}=m_{\alpha}=f$. However, since $\deg(f)\le \deg(F\diamond f)$, it follows that $$k=\deg(f)\le \deg(f_1)\le \ldots \le \deg(f_m)=k,$$ and so $\deg(f_1)=\deg(F\diamond f)=k$. Combining this last argument with Proposition~\ref{prop:degree-preserving} and Lemma~\ref{lem:trivials}, we obtain the following theorem.

\begin{theorem}\label{thm:degree-preserve}
Let $\I$ be the set of all monic irreducible polynomials over $\F_q$ and $F(x)\in \F_q[x]$ is a polynomial of degree $\ge 1$. Then the induced map $f\mapsto F\diamond f$ on  $\I$ preserves the degree of any irreducible polynomial in $\I$ if and only if $F(x)=ax^{p^h}+b$ for some $a, \in \F_q^*$,  $b\in \F_q$, and $h\ge 0$. 
\end{theorem}

Proposition~\ref{prop:degree-preserving} implies that if the map $f\mapsto F\diamond f$ on $\I$  is degree preserving then $F$ permutes the field $\overline{\F}_q$.  In the following, we study maps induced by permutations of $\F_{q^k}$, which are not neccessarily permutations of $\overline{\F}_q$.

\begin{theorem}\label{thm:local}
Let $k$ be a positive integer and  $\I_k$ be the set of monic irreducible polynomials over $\F_q$ of degree $k$.  Let $F(x)\in \F_q[x]$ such that the evaluation map $c\mapsto F(c)$ of $F$ on $\F_{q^k}$ is a permutation. Then for any $f\in \I_k$,  the polynomial $F\diamond f$ is also in $\I_k$, i.e., the restriction of the map $f\mapsto F\diamond f$ to the set $\I_k$ is a degree preserving map. Moreover, this restriction is also a permutation of the set $\mathcal{I}_k$.
\end{theorem}

We observe that, since $\I_k$ comprises the minimal polynomials of the elements in $\C_k$, it is sufficient to prove that $F$ permutes the set $\C_k$: in fact, if this occurs, $f\mapsto F\diamond f$ maps the set $\I_k$ to itself. In addition, if $f, g\in \I_k$ and $\alpha, \beta\in \C_k$ are roots of $f$ and $g$, respectively, then $F\diamond f=F\diamond g$ implies that $F(\alpha)$ and $F(\beta)$ have the same minimal polynomial. In particular, $F(\alpha)=F(\beta)^{q^i}=F(\beta^{q^i})$ for some $i\ge 0$ and, since $F$ is a permutation of $\C_k$, it follows that $\alpha=\beta^{q^i}$ and so $\alpha$ and $\beta$ are conjugates. Therefore, $\alpha$ and $\beta$ have the same minimal polynomial over $\F_q$, i.e., $f=g$. In other words, $f\mapsto F\diamond f$ maps $\I_k$ (a finite set) into itself and is an one-to-one map, hence is a permutation.

Next we show that $F$ permutes the set $\C_k$  and so we finish the  proof of Theorem~\ref{thm:local}.

\begin{proposition}\label{local-permut}
Let $F\in \F_q[x]$ and let $k$ be a positive integer such that the evaluation map $c\mapsto F(c)$ of $F$ on $\F_{q^k}$ is a permutation. Then the restriction of this evaluation map on the set $\C_k$ is a permutation of $\C_k$.
\end{proposition}

\begin{proof}
Let $p_1, \ldots, p_s$ are the distinct prime divisors of $k$, hence $$\bigcup_{i=1}^s \F_{q^{k/p_i}}=\F_{q^k}\setminus \C_k.$$ We observe that, for any positive integer $d$, $F(\F_{q^d})\subseteq \F_{q^d}$, hence $$\bigcup_{i=1}^s \F_{q^{k/p_i}}\supseteq F\left(\bigcup_{i=1}^s \F_{q^{k/p_i}}\right),$$ and, since $F$ is a permutation of $\F_{q^k}$, the previous inclusion is an equality of sets and then $F(\C_k)=\C_k$, i.e., $F$ permutes the set $\C_k$.
\end{proof}

\section{Permutations of irreducible polynomials}
Theorem~\ref{thm:local} says that, for a permutation polynomial $F\in \F_q[x]$ of $\F_{q^k}$, the map $f\mapsto F\diamond f$ is a permutation  of the set $\I_k$. In this context, it is interesting to study the permutations of $\F_{q^k}$ that are induced by polynomials in $\F_q[x]$. We first observe that two polynomials $F$ and $F_0$ induce the same evaluation map in $\F_{q^k}$ if and only if $F\equiv F_0\pmod{x^{q^k}-x}$. We consider the following set
$$\G_k:=\{P\in \F_q[x]\,|\, P\; \text{is a permutation polynomial of}\,\, \F_{q^k}\, \text{and}\, \deg(P)<q^k\}.$$
We shall prove that the set $\G_k$ has a group structure. First, we introduce a simple (but powerful) result.
\begin{proposition}[Frobenius-Stable Polynomials]\label{prop:frobenius} Let $f\in \overline{\F}_q[x]$ be a polynomial of degree $n$. The following are equivalent.
\begin{enumerate}[(i)]
\item There exists a set $C\subseteq \overline{\F}_q$ with cardinality at least $n+1$ such that $$f(\alpha^q)=f(\alpha)^q, \alpha\in C.$$
\item The coefficients of $f$ lie in $\F_q$, i.e., $f\in \F_q[x]$.
\end{enumerate}
\end{proposition}
\begin{proof}
We observe that, if $f(x) \in \F_q[x]$, then $f(\alpha^q)=f(\alpha)^q$ for any $\alpha\in \overline{\F}_q$ so it suffices to prove that (i) implies (ii). For this, let $C$ be as above and write $f(x)=\sum_{i=0}^{n}a_ix^i$. Consider the polynomial $f^*(x)=\sum_{i=0}^{n}(a_i^q-a_i)x^i$. Since $f(\alpha^q)=f(\alpha)^q$ for any $\alpha\in C$, we see that $f^*$ vanishes at every element of $C$. However, since the degree of $f^*$ is at most $n$ and $C$ has at least $n+1$ elements, it follows that $f^* = 0$, i.e., $a_i=a_i^q$. In other words, $f$ is a polynomial with coefficients in $\F_q$.
\end{proof}

As a consequence of the previous proposition, we obtain the following result.

\begin{corollary}\label{cor:group}
The set $\G_k$ is a group with respect to the composition modulo $x^{q^k}-x$.
\end{corollary}

\begin{proof}
Let $P, Q\in \G_k$. In particular, $P, Q\in \F_q[x]$ and so the reduction of $P\circ Q$ modulo $x^{q^k}-x$ is also a polynomial with coefficients in $\F_q$ and induces the same permutation of $P\circ Q$ in $\F_{q^k}$, i.e., $P\circ Q\in \G_k$. The identity element is the trivial permutation $P(x)=x\in \F_{q}[x]$. Observe that, from definition, any element of $\G_k$ has inverse (not necessarily in $\G_k$). Let $P\in \G_k$ and let $P_0$ be its inverse;  without loss of generality, $n=\deg P_0<q^k$. We just need to show that $P_0$ is a polynomial with coefficients in $\F_q$. Indeed, for any $\alpha\in \F_{q^k}$, we have $P(P_0(\alpha)^q) = P(P_0(\alpha))^q = \alpha^q = P(P_0(\alpha^q))$ and so $P_0(\alpha)^q = P_0(\alpha^q)$ since $P \in \G_k$. Since $|\F_{q^k}|=q^k\ge n+1$, from Proposition~\ref{prop:frobenius}, it follows that $P_0\in \F_q[x]$.
\end{proof}

If $\sym(\I_k)$ denotes the symmetric group of the set $\I_k$, from Theorem~\ref{thm:local}, we have the group homomorphism $\tau_k:\G_k\to \sym(\I_k)$ given as follows: for $P\in \G_k$, $\tau_k(P)=\varphi_P$, where $\varphi_P:\I_k\to \I_k$ is such that $\varphi_P(f)=P\diamond f$. The following theorem shows that this homomorphism is \emph{onto} and, in particular, this implies that any permutation of the set $\I_k$ can be viewed as a map $f\mapsto P\diamond f$ for some $P\in \G_k$.

\begin{theorem}\label{thm:onto}
For any permutation $\sigma$ of the set $\I_k$, there exists an element $P\in \G_k$ such that $\sigma(f)=P\diamond f$ for any $f\in \I_k$.
\end{theorem}

\begin{proof}
Let $f_1, \ldots, f_{n_k}$ be a list of all elements in the set $\I_k$, where $n_k=|\I_k|$. For each $1\le i\le n_k$, let $\alpha_i$ be any root of $f_i$. We observe that $\C_k$ comprises the elements $\alpha_i^{q^j}$ with $1\le i\le n_k$ and $0\le j\le k-1$. Fix $\sigma\in \sym(\I_k)$. Then $\sigma$ induces a permutation $\lambda$ of the set $\{1, \ldots, n_k\}$ such that $\sigma(f_i)=f_{\lambda(i)}$. Let $P\in \F_{q^k}[x]$ be the polynomial of least degree such that 
$$P\left(\alpha_i^{q^j}\right)=\alpha_{\lambda(i)}^{q^j}, \,\text{for any}\; 1\le i\le n_k\; \text{and} \;0\le j\le k-1,$$
and $P(\beta)=\beta$ if $\beta\in \F_{q^k}\setminus \C_k$. Using \emph{Lagrange interpolation}, such a $P$ exists and has degree at most $q^k-1$. From definition, $P$ is a permutation of $\F_{q^k}$ and we can verify that $P(\alpha^q)=P(\alpha)^q$ for any $\alpha\in \F_{q^k}$. From Proposition~\ref{prop:frobenius}, it follows that $P\in \F_q[x]$. In conclusion, $P$ is an element of $\G_k$. We observe that, for each $1\le i\le n_k$, $P(\alpha_i)=\alpha_{\lambda(i)}$ and so
$$P\diamond f_i=m_{\alpha_{\lambda(i)}}(x)=f_{\lambda(i)}=\sigma(f_i).$$
\end{proof}

We observe that, for polynomials $F\in \F_q[x]$ and $f\in \I_k$, $F\diamond f$ is the minimal polynomial of $F(\alpha)$, where $\alpha$ is any root of $f$. In particular, the computation of $F\diamond f$ requires the construction of the field $\F_{q^k}$.  
When $F=P$ is a permutation of $\F_{q^k}$, the following lemma gives an alternative way of obtaining $P\diamond f$.

\begin{lemma}\label{lem:gcd}
Let $P\in \G_k$ and let $Q\in \G_k$ be the inverse of $P$. For any $f\in \I_k$, the polynomial $f(Q(x))$ is such that any of its irreducible factors over $\F_q$ has degree divisible by $k$. Additionally, $P\diamond f$ is the only irreducible factor of degree $k$ of $f(Q(x))$ (possibly with multiplicity greater than one). In particular,
\begin{equation}\label{eq:gcd-action-0}P\diamond f(x)=\gcd(f(Q(x)), x^{q^k}-x).\end{equation}
\end{lemma}

\begin{proof}
Let $g$ be any irreducible factor of $f(Q(x))$ and let $\beta\in \overline{\F}_q$ be any element such that $g(\beta)=0$. We observe that $f(Q(\beta))=0$ and so, there exists a root $\alpha\in \F_{q^k}$ of $f$ such that $Q(\beta)=\alpha$. The last equality says that $\beta$ is a root of $Q(x)-\alpha$. Since $\alpha$ is an element of degree $k$, we conclude that $\beta$ is in an extension of $\F_{q^{k}}$, hence $\deg(\beta)=kd$ for some $d\ge 1$. Clearly $g$ is the minimal polynomial of $\beta$ over $\F_q$ and so $\deg(g)=kd$. In particular, if $\deg(g)=k$, then $\beta\in \F_{q^k}$ and, since $Q\in \G_k$ is a permutation of $\F_{q^k}$ with $Q(\beta)=\alpha$, it follows that $\beta=P(\alpha)$ and, from definition, $g=P\diamond f$. This shows that $P\diamond f$ is the only factor of degree $k$ of $f(Q(x))$. Since any other factor of $f(Q(x))$ has degree $kd$ for some $d>1$ and $x^{q^k}-x$ has no repeated irreducible factors, it follows that $P\diamond f=\gcd(f(Q(x)), x^{q^k}-x)$.
\end{proof}

Though the computation of GCD's in $\F_{q}[x]$ does not require the construction of the field $\F_{q^k}$, the explicit computation of inverse of a permutation is, in general, a hard problem. For this reason, we introduce the following alternative operation.

\begin{define}
For an element $P\in \G_k$ and $f\in \I_k$, we set 
$$P\ast f=Q\diamond f,$$
where $Q\in \G_k$ is the compositional inverse of $P$.
\end{define}

One can easily see that, if $\alpha$ is any root of $f\in \I_k$, $P\ast f=m_{\beta}$, where $\beta=Q(\alpha)$ is the only element in $\F_{q^k}$ such that $P(\beta)=\alpha$. For an element $P\in \G_k$, the compositions $P\ast f$ and $P\diamond f$ are \emph{dual}, in the sense that $f=P\diamond (P\ast f)=P\ast (P\diamond f)$. The advantage is that the composition $P\ast f$ can be easily computed; from Lemma~\ref{lem:gcd}, it follows that
\begin{equation}\label{eq:gcd-action}P\ast f(x)=\gcd(f(P(x)), x^{q^k}-x),\end{equation}
for any $P\in \G_k$ and $f\in \I_k$.

\begin{remark}
We emphasize that there is no loss of generality on working with the compositions $P\ast f$: in fact, the dynamics of $P$ on $\I_k$ via the compositions $P\ast f$ and $P\diamond f$ are essentially the same.
\end{remark}

\noindent For the rest of this paper, we consider the maps $f\mapsto P\ast f$ induced by elements $P\in \G_k$ on the set $\I_k$.

\begin{example}\label{ex:1}
Let $q=2$ and $k=4$. It is direct to verify that $P=x^7$ is a permutation of $\F_{2^4}=\F_{16}$.  We have $\I_4=\{f_1, f_2, f_3\}$, where $f_1(x)=x^4+x+1, f_2(x)=x^4+x^3+1$ and $f_3(x)=x^4+x^3+x^2+x+1$.  Using the formula given in Eq.~\eqref{eq:gcd-action}, we obtain $P\ast f_1=f_2$, $P\ast f_2=f_1$ and $P\ast f_3=f_3$. This corresponds to the permutation of three symbols with cycle decomposition $(1\, 2)\, (3)$.

\end{example}

\begin{example}\label{ex:2}
Let $q=3$ and $k=3$. It is direct to verify that, for $h=x+1$, $L_h=x^3+x$ is a permutation of $\F_{3^3}=\F_{27}$.  We have $\I_3=\{f_1, \ldots, f_8\}$, where $f_1(x)=x^3-x+1, f_2(x)=x^3-x-1, f_3(x)=x^3+x^2-1, f_4(x)=x^3+x^2+x-1, f_5(x)=x^3+x^2-x+1, f_6(x)=x^3-x^2+1, f_7(x)=x^3-x^2+x+1$ and $f_8(x)=x^3-x^2-x-1$.  Using the formula given in Eq.~\eqref{eq:gcd-action}, we obtain $P\ast f_1=f_2$, $P\ast f_2=f_1$,  $P\ast f_3=f_8$, $P\ast f_8=f_4$, $P\ast f_4=f_6$, $P\ast f_6=f_5$, $P\ast f_5=f_7$ and $P\ast f_7=f_3$. This corresponds to the permutation of eight symbols with cycle decomposition $(1\, 2)\, (3\, 8\, 4\, 6\, 5\, 7)$.
\end{example}

\subsection{M\"{o}bius maps on irreducible polynomials}\label{subsec:Mobius-action}
If $\gamma(x)=\frac{ax+b}{cx+d}$ is a M\"{o}bius map with $a, b, c, d\in \F_q$ and $ad-bc\ne 0$, we have that $\gamma$ is a permutation of $\F_{q^k}$ for any $k$ (with a suitable evaluation at the possible pole of $\gamma(x)$): its inverse is given by $\gamma^{-1}(x)=\frac{dx-b}{a-cx}$. If $f\in \I_k$ with $k\ge 2$ and $\alpha\in \F_{q^k}\setminus \F_q$ is any root of $f$, the minimal polynomial of $\gamma^{-1}(\alpha)=\frac{d\alpha-b}{a-c\alpha}$ over $\F_q$ equals $c_{f}(cx+d)^{k}f\left(\frac{ax+b}{cx+d}\right)$, where $c_f$ is the only nonzero element of $\F_{q}$ such that $c_{f}(cx+d)^{k}f\left(\frac{ax+b}{cx+d}\right)$ is a monic polynomial. In particular, \begin{equation}\label{eq:mobius-map}\gamma \ast f=c_{f}(cx+d)^{k}f\left(\frac{ax+b}{cx+d}\right).\end{equation}
M\"{o}bius transformations on irreducible polynomials of degree $k\ge 2$ like the one given in Eq.~\eqref{eq:mobius-map} are considered in \cite{ST12}.

\section{Invariant theory for the maps $f\mapsto P\ast f$}
 In this section, we provide a general invariant theory on the dynamics of $P\in \G_k$ on $\I_k$ via the map $f\mapsto P\ast f$.  The study of \emph{fixed points} plays a main role in the theory of dynamics. When considering dynamics on finite sets, the number of fixed points is frequently considered. Throughout this section, $k$ is a fixed positive integer.

\begin{define}Given $P\in \G_k$, $\mathcal C_P=\{f\in \I_k\,|\, P\ast f=f\},$ is the set of fixed points of $\I_k$ by $P$ and $n_P=|\mathcal C_P|$ is the number of fixed points.\end{define}

In the following theorem, we give a simple characterization of polynomials $f\in \I_k$ that are fixed by an element $P\in \G_k$.

\begin{theorem}\label{thm:1}
For $f\in \I_k$ and $P\in \G_k$, the following are equivalent:
\begin{enumerate}[(i)]
\item $P\ast f=f$, i.e., $f\in \mathcal C_P$;
\item $f(x)$ divides $x^{q^i}-P(x)$ for some $0\le i\le k-1$.
\end{enumerate}
In particular, if we set $R_d[P](x)=\prod_{i=0}^{d-1}(x^{q^i}-P(x))$, then $P\ast f=f$ if and only if $f$ divides $R_k[P]$.
\end{theorem}

\begin{proof}
Let $\alpha\in \F_{q^k}$ be any root of $f$ and let $\beta\in \F_{q^k}$ be the element such that $P(\beta)=\alpha$. We observe that $P\ast f=f$ if and only if the minimal polynomial of $\beta$ is $f$, i.e., $f(\beta)=0$. In other words, $P\ast f=f$ if and only if $\beta=\alpha^{q^j}$ for some $0\le j\le k-1$. The latter holds if and only if $P(\alpha^{q^j})=\alpha$, i.e., $f$ divides $x^{q^{k-j}}-P(x)$, where $0\le k-j\le k-1$.
\end{proof}

\begin{define}
For each positive integer $d$, we set $$\Psi_d(x)=\prod_{f\in \I_d}f(x)=\prod_{\alpha\in \C_d}(x-\alpha).$$ 
\end{define}

It is clear that $x^{q^k}-x=\prod_{d|k}\Psi_d(x)$. From the previous theorem, a general implicit formula for the number of fixed points can be derived.

\begin{theorem}
For any polynomial $P$ and integers $i\ge 0$,  $d\ge 1$, set $$R_d[P](x)=\prod_{i=0}^{d-1}(x^{q^i}-P(x)), \; g_P^{(i, d)}(x)=\gcd(x^{q^i}-P(x), x^{q^d}-x),$$ and $h_P^{(d)}(x)=\gcd(R_d[P](x), \Psi_d(x))$. If $P\in \G_k$, the number $n_P$ of fixed points of $\I_k$ by $P$ satisfies the following identity

\begin{equation}\label{fix-points}n_P=\frac{1}{k}\deg(h_P^{(k)})=\sum_{d|k}\frac{\mu(k/d)}{d}\sum_{i=0}^{d-1}\deg(g_P^{(i, d)}).\end{equation}
\end{theorem}

\begin{proof}
Since $\Psi_k$ is squarefree, the equality $n_P=\frac{1}{k}\deg(h_P^{(k)})$ follows directly from Theorem~\ref{thm:1}. Notice that any irreducible factor of $g_P^{(i, k)}$ is of degree a divisor of $k$. Let $F$ be any monic irreducible polynomial of degree $d$, where $d$ is a divisor of $k$. We observe that $F$ divides $g_P^{(i, k)}$ if and only if $F$ divides $g_P^{(i, d)}$. In this case, the minimal $i_0$ with such property satisfies $i_0\le d-1$, i.e., $F$ divides $R_d[P]$. Additionally, if such $i_0$ exists, $F$ divides $g_p^{(j, k)}$ with $0\le j\le k-1$ if and only if $j=i_0+ds$ with $0\le s< \frac{k}{d}$.  Indeed, $F$ divides $x^{q^{j}-q^{i_0}}-1$ if and only if $F$ divides $x^{q^{j-i_0}-1}-1$, or $x^{q^{j-i_0}}-x$. Because $F$ is irreducible and has degree $d$, we must have $j\equiv i_0\pmod d$. In particular, if $F$ is an irreducible polynomial of degree $d$ and divides $R_k[P]$, then  $F^{k/d}$ is the greatest power of $F$ that divides $R_k[P]$; this implies that
$$\prod_{i=0}^{k-1}\gcd(x^{q^i}-P(x), \Psi_d(x))=\gcd\left(\prod_{i=0}^{d-1}(x^{q^i}-P(x)), \Psi_d(x)\right)^{k/d}=(h_P^{(d)}(x))^{k/d}.$$
Therefore, since $x^{q^k}-x=\prod_{d|k}\Psi_d(x)$ is squarefree, the following holds:

\begin{equation}\label{mobius}\prod_{i=0}^{k-1}g_P^{(i, k)}(x)=\prod_{i=0}^{k-1}\prod_{d|k}\gcd(x^{q^i}-P(x), \Psi_d(x))=\prod_{d|k}(h_P^{(d)}(x))^{k/d}.\end{equation}
We observe that Eq.~\eqref{mobius} holds for any positive integer $k$ and any polynomial $P$. For a positive integer $s$, we set $\mathcal L(s)=\frac{1}{s}\sum_{i=0}^{s-1}\deg(g_P^{(i, s)})$ and $\mathcal M(s)=\frac{1}{s}\deg (h_P^{(s)})$. Taking degrees on Eq.~\eqref{mobius}, we see that $\mathcal L(k)=\sum_{d|k}\mathcal M(d)$ for any positive integer $k$. From the \emph{M\"{o}bius inversion formula}, it follows that

$$n_P=\frac{1}{k}\deg(h_P^{(k)})=\mathcal M(k)=\sum_{d|k}\mathcal L(d)\cdot \mu(k/d).$$
\end{proof}

When the permutation $P$ is a monomial or a linearized polynomial, Eq.~\eqref{fix-points} can be fairly simplified.

\begin{corollary}\label{cor:fix-points-lin-monomial}
The following hold:

\begin{enumerate}[(i)]
\item If $P(x)=x^n$ is a permutation polynomial of $\F_{q^k}$, i.e., $\gcd(n, q^k-1)=1$, then
\begin{equation}\label{eq:monomial-fix}n_P=\varepsilon(k)+\sum_{d|k}\frac{\mu(k/d)}{d}\sum_{i=0}^{d-1}\gcd(q^i-n, q^d-1),\end{equation}
where $\varepsilon(k)=0$ if $k\ne 1$ and $\varepsilon(1)=1$.

\item Let $h(x)\in \F_q[x]$ with $h(x)=\sum_{i=0}^{k-1}a_ix^i$ be a polynomial such that $L_h(x)=\sum_{i=0}^{k-1}a_ix^{q^i}$ is a permutation of $\F_{q^k}$, i.e., $\gcd(h(x), x^k-1)=1$. For $P=L_h$, we have
\begin{equation}\label{eq:linearized-fix}n_P=\sum_{d|k}\frac{\mu(k/d)}{d}\sum_{i=0}^{d-1}q^{r_{i, d}},\end{equation}
where $r_{i, d}=\deg(\gcd(x^i-h, x^d-1))$.
\end{enumerate}
\end{corollary}

\begin{proof}
Applying Eq.~\eqref{fix-points}, we just need to compute $\deg(g_P^{(i, d)})$ explicitly.
\begin{enumerate}[(i)]
\item We observe that, for $P=x^n$, the following holds: 
$$g_P^{(i, d)}(x)=\gcd(x^{q^i}-x^n, x^{q^d}-x)=x\cdot \gcd(x^{q^i-n}-1, x^{q^d}-1)=x\cdot (x^{\gcd(q^i-n, q^d-1)}-1).$$ 
Therefore, $\deg(g_P^{(i, d)})=1+\gcd(q^i-n, q^d-1)$. To finish the proof, we observe that $\sum_{d|k}\mu(\frac{k}{d})=\varepsilon(k)$ for any positive integer $k$.

\item From Lemma~\ref{lem:linearized-properties}, if $f$ and $g$ are polynomials with coefficients in $\F_q$, then $\gcd(L_f, L_g)=L_{\gcd(f, g)}$. In particular, for $P=L_h$, we have $$g_P^{(i, d)}(x)=\gcd(x^{q^i}-L_h, x^{q^d}-x)=\gcd(L_{x^i-h}, L_{x^d-1})=L_{\gcd(x^i-h, x^d-1)},$$ and so 
$\deg(g_P^{(i, d)})=q^{r_{i, d}}$, where $r_{i, d}=\deg(\gcd(x^i-h, x^d-1))$.
\end{enumerate}
\end{proof}

\begin{example}
We consider the monomial permutation polynomial $P(x)=x^7$ over $\F_{2^4}=\F_{16}$. From the previous corollary, the number of fixed points of $P$ is

$$n_P=\sum_{d|4}\frac{\mu(4/d)}{d}\sum_{i=0}^{d-1}(\gcd(2^i-7, 2^d-1)+1)=1,$$
as expected by Example~\ref{ex:1}.
\end{example}

In the following proposition we show that, in some special extensions of finite fields, the number of fixed points can be given explicitly when considering monomial and linearized permutations.

\begin{proposition}
The following hold:

\begin{enumerate}[(i)]
\item If $k$ and $r:= \frac{q^k-1}{q-1}$ are prime numbers and $P(x)=x^n\in \G_k$ is a permutation polynomial of $\F_{q^k}$, then the number of fixed points of $\I_k$ by $P$ is 

\begin{eqnarray*}
n_P = \left\{
\begin{array}{ll}
\frac{(r-1)}{k} \gcd(n-1, q-1), &  \text{if~}  n\equiv q^i ~(\bmod~{r}) \text{~for some~} 0 \leq i \leq k-1; \\
0, & \text{otherwise}.
\end{array}
\right.
\end{eqnarray*}

\item If $k$ is a prime number, $T(x):=\frac{x^k-1}{x-1}=x^{k-1}+\cdots+x+1$ is an irreducible polynomial and $P(x)=L_f(x) \in \G_k$ is a permutation polynomial of $\F_{q^k}$ with $f(x)\in \F_q[x]$, then we have 

{\small

\begin{eqnarray*}
n_P = \left\{
\begin{array}{ll}
\frac{q^2-q}{2}=|\I_2|, &  \text{if~}  q  \text{~is even}, k = 2 \text{ and } f(x) = 1 \text{ or } x;  \\
\frac{q^k-q}{k}=|\I_k|,   &  \text{if~} q \text{~is odd or } k \neq 2\text{ and } f(x) = x^i \text{~for ~} 0 \leq i \leq k-1; \\
\frac{q^{k-1}-1}{k},  &   \text{if~} q \text{~is odd or } k \neq 2,    f(x) = a T(x) + x^i \text{~for ~} 0 \leq i \leq k-1, a \in \fq^*;  \\

0, & otherwise.
\end{array}
\right.
\end{eqnarray*}

}



\end{enumerate}
\end{proposition}

\begin{proof}
\begin{enumerate}[(i)]
\item Since $k\ge 2$ and $r$ is prime, we have $r>q+1$ and so $1=\gcd(r, q-1)=\gcd(k, q-1)$. Since $P=x^n\in \G_k$ is a monomial permutation, $\gcd(n, q^k-1)= 1$ and $n<q^k$. In particular, $n=ar+b$ with $a\le q-1$ and $b<r$. Set $s=\gcd(n-1, q-1)$ and, for each positive integer $0\le i\le k-1$, set $s_i=\gcd(q^i-n, r)$. Therefore, $\gcd(q^i-n, q^k-1)=\gcd(q^i-n, r)\cdot \gcd(q^i-n, q-1)=s_i\cdot \gcd(n-1, q-1)=s_is$. Since $k$ is prime, from Eq.~\eqref{eq:monomial-fix}, we obtain the following equality:
$$n_P=-s+\frac{1}{k}\sum_{i=0}^{k-1}s_is=\frac{s}{k}\sum_{i=0}^{k-1}(s_i-1).$$
We first observe that if $N$ is a positive integer that divides $s_i$ and $s_j$, then $N$ divides $q^i-q^j$, where $|q^i-q^j|\le q^{k-1}-1<r$. In particular, for at most one index $0\le i\le k-1$, we have $s_i\ne 1$. If $s_i=1$ for every $0\le i\le k-1$, the previous equality yields $n_P=0$.
If there exists $0\le i\le k-1$ such that $s_i\ne 1$, then $q^i-n$ is divisible by $r$. However, since $n=ar+b$, it follows that $q^i-b$ is divisible by $r$. Because $0\le q^i, b<r$, the latter implies $q^i-b=0$. Hence $n=ar+q^i=a\cdot \frac{q^k-1}{q-1}+q^i$ and, in this case, $s_j=1$ if $j\ne i$ and $s_i=r$. In addition, $s=\gcd(n-1, q-1)=\gcd(ar+q^i-1, q-1)=\gcd(ar, q-1)=\gcd(a, q-1)$, since $\gcd(r, q-1)=1$. Therefore, $n_P=\frac{s}{k}(r-1)=\gcd(a, q-1)\cdot \frac{q^k-q}{k(q-1)}$.

\item We split the proof into cases.

\begin{itemize}\item If $q$ is even and $k=2$, $x-1=x+1$ and $x^2-1=(x+1)^2$. Since $P=L_f\in \G_k=\G_2$, $\deg(L_f)<q^2$ and so $f$ is a linear polynomial. If $G=\gcd(f(x)-1, x+1), G_1=\gcd(f(x)-1, (x+1)^2)$ and $G_2=\gcd(f(x)-x, (x+1)^2)$ are polynomials of degrees $m, m_1$ and $m_2$, respectively, from Eq.~\eqref{eq:linearized-fix}, we obtain the following equality:
$$n_P=-q^m+\frac{q^{m_1}+q^{m_2}}{2}.$$ Since $f(x)-x\equiv f(x)-1\pmod {x+1}$, one can see that the polynomials $G, G_1$ and $G_2$ are all equal to $1$ or all distinct from $1$. In particular, the numbers $m, m_1$ and $m_2$ are all zero or all nonzero. If they are all zero, $c_P=0$. Suppose now that $m, m_1$ and $m_2$ are not zero. Since $f$ is a linear polynomial, $m, m_1$ and $m_2$ are at most one unless $f=1$ or $f=x$ and, in this case, $m=1$, and $\{m_1, m_2\}=\{1, 2\}$ (in some order). For $m=m_1=m_2=1$, we obtain $n_P=0$. Otherwise, $n_P=-q+\frac{q^2+q}{2}=\frac{q^2-q}{2}=|\I_2|$.

\item  If $q$ is odd or $k\ne 2$, we observe that $T(x)=\frac{x^k-1}{x-1}$ is divisible by $x-1$ if and only if $T(1)=k$ equals zero, i.e., $k$ is divisible by $p$, the characteristic of $\F_{q}$. Since $k$ is prime, $k=p$, hence $T(x)=(x-1)^{p-1}$ which is not irreducible if $p>2$, a contradiction. Therefore $T(x)$ and $x-1$ are relatively prime.

Since $P=L_f\in \G_k$, $\deg(L_f)<q^k$ and thus $f$ is a polynomial of degree at most $k-1$. Write $f=bT+R$, where $b\in \F_q$ and $R$ is a polynomial of degree at most $k-2$.  Set $G=\gcd(f-1, x-1)$ and, for each positive integer $0\le i\le k-1$, set $G_i=\gcd(x^i-f, T)$. Therefore, $\gcd(x^i-f, x^k-1)=\gcd(x^i-f, T)\cdot \gcd(x^i-f, x-1)=G_i\cdot \gcd(f-1, x-1)=G_iG$. If we set $m_i=\deg(G_i)$, from Eq.~\eqref{eq:linearized-fix}, we obtain the following equality:
$$n_p=-q^{m}+\frac{1}{k}\sum_{i=0}^{k-1}q^{m_i+m}=\frac{q^m}{k}\sum_{i=0}^{k-1}(q^{m_i}-1).$$
We observe that if $A(x)$ is a polynomial that divides $G_i$ and $G_j$, then $A(x)$ divides $x^i-x^j$, which is a polynomial of degree at most $k-1$ (and distinct from $T(x)$ since $q$ is odd or $k\ne 2$). In particular, for at most one index $0\le i\le k-1$, we have $G_i\ne 1$, i.e., $m_i\ne 0$. If $G_i=1$ for every $0\le i\le k-1$, then $m_i=0$ and the previous equality yields $n_P=0$. If there exists $0\le i\le k-1$ such that $G_i\ne 1$, then $x^i-f$ is divisible by $T$: since $f=bT+R$, it follows that $x^i-R$ is divisible by $T$, where $R$ has degree at most $k-2$. If $i\le k-2$, then $x^i-R$ is a polynomial of degree at most $k-2$ and the latter implies $x^i-R=0$, i.e., $f(x)=bT(x)+x^i$. If $i=k-1$, then it follows that $x^{k-1}-R=T$, hence $f=bT+R=(b-1)T+x^{k-1}$. In particular, if $n_P\ne 0$, then $f$ is of the form $aT+x^i$ for some $0\le i\le k-1$ and, in this case, $m_j=0$ if $j\ne i$ and $m_i=k-1$.Therefore, $$n_P=q^{m}\frac{1}{k}\sum_{i=0}^{k-1}(q^{m_i}-1)=\frac{q^{m}(q^{k-1} -1) }{k}.$$ Now we just need to compute $m$. If $f(x)=aT(x) +x^i$, then $f(1)=aT(1)+1=ak+1$. Hence $G=\gcd(f-1, x-1)=1$ unless $ak+1=1$, i.e., $ak=0$. We have seen that $k$ is not divisible by $p$, hence $a=0$. In conclusion, $m=\deg(G)$ is zero unless $a=0$ and, for $a=0$, we must have $m=1$. This concludes the proof.
\end{itemize}
\end{enumerate}
\end{proof}

\subsection{Invariants via  special permutations}
So far we have considered polynomials $f$ that are fixed by an element $P\in \G_k$, i.e., $P\ast f=f$. One may ask if there is any general characteristic of $P\ast f$ that is \emph{inherited} from $f$. In other words, we are interested in the algebraic polynomial structures that are preserved by the compositions $P\ast f$. As follows, we show that monomial and linearized permutations preserve certain algebraic structures of polynomials.

\subsubsection{Multiplicative properties invariant by monomial permutations}
Here we employ many definitions and notations introduced in Subsection~\ref{subsec:mult}.

\begin{proposition}
Suppose that $f\in \I_k$ with $f\ne x$ and $P=x^n$ is a monomial permutation over $\F_{q^k}$. Then $\ord(f)=\ord(P\ast f)$, i.e., the permutation of $\I_k$ induced by $P=x^n$ preserves the multiplicative order of polynomials. In particular, if $f$ is a primitive polynomial, then so is $P\ast f$.
\end{proposition}

\begin{proof}
Let $\alpha$ be any root of $f$ and let $\beta$ be the unique element of $\F_{q^k}$ such that $\beta^n=\alpha$. We know that $\gcd(n, q^k-1)=1$ and then, from item (ii) of Lemma~\ref{thm:mult}, it follows that $\ord(\beta)=\ord(\alpha)$. To finish the proof, just recall that $P\ast f=m_{\beta}(x)$, $\ord(f)=\ord(\alpha)$ and $\ord(m_{\beta})=\ord(\beta)$. 
\end{proof}

\begin{define} For an element $f\in \I_k$ such that $f(x)=x^k+\sum_{i=0}^{k-1}a_ix^i$, the element $a_0$ is the \emph{norm} of $f$. \end{define}
We can easily verify that, if $\alpha$ is any root of $f$, then $a_0=\prod_{0\le i\le k-1}\alpha^{q^i}=\alpha^{\frac{q^k-1}{q-1}}$. In the following proposition, we show how the norms of $f$ and $P\ast f$ are related, in the case when $P$ is a monomial.

\begin{proposition}
Suppose that $f\in \I_k$, $P=x^n$ is a monomial permutation of $\F_{q^k}$ (i.e., $\gcd(n, q^k-1)=1$) and let $n_0<q^k-1$ be the positive integer such that $nn_0\equiv 1\pmod {q^k-1}$. If $f$ has norm $a\in \F_q$, then $x^n\ast f$ has norm $a^{n_0}$. In particular, the following hold.
\begin{enumerate}[(i)]
\item if $f$ has norm $1$, then so has $P\ast f$;
\item if $n\equiv 1\pmod {q-1}$, the polynomials $f$ and $P \ast f$ have the same norm for every $f\in \I_k$. 
\end{enumerate}
\end{proposition}

\begin{proof}
Let $\alpha$ be any root of $f$ and let $\beta$ be the unique element of $\F_{q^k}$ such that $\beta^n=\alpha$, hence $\beta=\alpha^{n_0}$. Hence $f=m_{\alpha}$ and $x^n\ast f=m_{\beta}$. With this notation, $f$ has norm $a=\alpha^{\frac{q^k-1}{q-1}}$ and
$x^n\ast f$ has norm $\beta^{\frac{q^k-1}{q-1}}=a^{n_0}$.  From this fact, items (i) and (ii) are straightforward to check.
\end{proof}

In particular, the previous proposition entails that, if $P=x^n$ is a monomial permutation such that $n\equiv 1 \pmod {q-1}$, then $f$ and $P\ast f$ have the same constant coefficient, for any $f\in \I_k$.

\subsubsection{Additive properties invariant by linearized permutations}
Here we employ many definitions and notations introduced in Subsection~\ref{subsec:add}.
\begin{proposition}
Suppose that $f\in \I_k$ and $L_g$ is a permutation over $\F_{q^k}$, where $L_g$ is the $q$-associate of $g\in \F_q[x]$. Then $\oord(f)=\oord(L_g\ast f)$, i.e., the permutation of $\I_k$ induced by $L_g$ preserves the $\F_q$-order of polynomials. In particular, if $f$ is a normal polynomial, then so is $L_g\ast f$.
\end{proposition}

\begin{proof}
Let $\alpha$ be any root of $f(x)$ and let $\beta$ be the unique element of $\F_{q^k}$ such that $L_g(\beta)=\alpha$. We know that $\gcd(g, x^k-1)=1$ and then, from item (ii) of Theorem~\ref{thm:add}, it follows that $m_{\alpha, q}(x)=m_{\beta, q}(x)$. To finish the proof, we recall that $L_g\ast f=m_{\beta}(x)$, $\oord(f)=m_{\alpha, q}$ and $\oord(m_{\beta})=m_{\beta, q}$. 
\end{proof}

\begin{define} For an element $f\in \I_k$ such that $f(x)=x^k+\sum_{i=0}^{k-1}a_ix^i$, $a_{k-1}\in \F_q$ is the \emph{trace} of $f(x)$ . \end{define}

We can easily verify that, if $\alpha$ is any root of $f$, then $a_{k-1}=\sum_{0\le i\le k-1}\alpha^{q^i}=L_T(\alpha)$, where $L_T$ is the $q$-associate of $T=\frac{x^k-1}{x-1}$. In the following proposition, we show how the traces of $f$ and $P\ast f$ are related, in the case when $P$ is linearized.

\begin{proposition}
Suppose that $f\in \I_k$ and $L_g$ is a permutation of $\F_{q^k}$, where $L_g$ is the $q$-associate of $g=\sum_{i=0}^{k-1}a_ix^i\in \F_q[x]$. Set $g(1)=c\in \F_q$. If $f$ has trace $a\in \F_q$, then $L_g\ast f$ has trace $c^{-1}a$. In particular, the following hold:

\begin{enumerate}[(i)]
\item if $f$ has trace $0$, then so has $L_g\ast f$.
\item if $g(1)=1$, then the polynomials $f$ and $L_g\ast f$ have the same trace for every $f\in \I_k$.
\end{enumerate}
\end{proposition}

\begin{proof}
Let $\alpha$ be any root of $f(x)$ and let $\beta$ be the unique element of $\F_{q^k}$ such that $L_g(\beta)=\alpha$. Hence $f=m_{\alpha}$ and $L_g\ast f=m_{\beta}$. With this notation, $f$ and $L_g\ast f$ have traces $L_T(L_g(\beta))$ and $L_T(\beta)$, respectively, where $T(x) =\frac{x^k-1}{x-1}$. If we set $L_T(\beta)=b\in \F_q$, we have that $b^{q^i}=b$ for any $i\ge 0$ and then $a=L_T(L_g(\beta))=L_g(L_T(\beta))=L_g(b)=\sum_{i=0}^{n-1}a_i b=bg(1)=bc$. Since $L_g$ is a permutation of $\F_{q^k}$, $\gcd(g, x^k-1)=1$ and, in particular, $g$ is not divisible by $x-1$ and so $g(1)=c\ne 0$. Therefore, $b=c^{-1}a$. From this fact, items (i) and (ii) are straightforward to check.
\end{proof}

In particular, the previous proposition says that, if $P=L_g$ is a linearized permutation such that $g(1)=1$, then $f$ and $P\ast f$ have the same coefficient of $x^{k-1}$, for any $f\in \I_k$.

\section{Dynamics of $\G_k$ on the sets $\C_k$ and $\I_k$}
If $F:S\to S$ is any map from a set $S$ to itself, we can associate to it a directed graph with nodes $\{a;\, a\in S\}$ and directed edges $\{(a, F(a);\, a\in S\}$: such a graph is called the {\it functional graph} of $F$ over $S$. The functional graph gives many informations of the dynamics of the function on the set: for instance, the \emph{orbit} $\{a, f(a), f(f(a)), \ldots\}$ of a point $a\in S$ is described by a \emph{path} in the functional graph. We know that any element $P\in \G_k$ induces permutations on the sets $\C_k$ and $\I_k$: namely, $P$ induces the evaluation map $c\to P(c)$ on $\C_k$ (see Proposition~\ref{local-permut}) and the map $f\mapsto P\ast f$ on $\I_k$. 

\begin{define}
For $P\in \G_k$, $G(P, \C_k)$ and $G(P, \I_k)$ denote the functional graphs of the evaluation map of $P$ on $\C_k$ and the map $f\mapsto P\ast f$ on $\I_k$, respectively.
\end{define}

We observe that $\I_k$ can be describe as the set of the minimal polynomials of the elements in $\C_k$. In this point of view, the following lemma shows how the graphs $G(P, \C_k)$ and $G(P, \I_k)$ interact. Its proof is straightforward so we omit.

\begin{lemma}
For any $\alpha, \beta \in \C_k$ and any $P\in \G_k$, the following hold:
\begin{enumerate}[(i)]
\item if $(\beta, \alpha)\in G(P, \C_k)$, then $(m_{\alpha}, m_{\beta})\in G(P, \I_k)$,
\item $(m_{\alpha}, m_{\beta})\in G(P, \I_k)$ if and only if $(\beta^{q^i}, \alpha)\in G(P, \C_k)$ for some $i\ge 0$ (or, equivalently, for some $0\le i\le k-1$).
\end{enumerate}
\end{lemma}

Since the map induced by a permutation $P\in \G_k$ on the set $\C_k$ (resp. $\I_k$) is one-to-one, any $a\in G(P, \C_k)$ (resp. any $f\in G(P, \I_k)$) is periodic. 

\begin{define}\label{def:periods}
Let $P\in \G_k$, $\alpha\in \C_k$ and $f\in \I_k$. 

\begin{enumerate}[(i)]
\item $c_P(\alpha)$ is the least period of $\alpha$ by $P$, i.e., $c_P(\alpha):=\min\{n>0\,|\, P^{(n)}(\alpha)=\alpha\}$. 
\item $c_P^*(f)$ is the least period of $f$ by $P$, i.e., $c_P^{*}(f):=\min \{n>0\,|\, P^{(n)}\ast f=f\}$. 
\item $S_P$ (resp. $S_P^*$) is the spectrum of the \emph{distinct} period lengths $c_P(\alpha), \alpha\in \C_k$ (resp. $c_P^*(f), f\in \I_k$).  
\item $\mu_k(P)=\min\{j\,|\, j\in S_P\}$ and $\mu_k^*(P)=\min\{j\,|\, j\in S_P^*\}$ are the minimal period lengths, respectively.
\end{enumerate}
\end{define}

In the following theorem, we present several relations between the numbers $c_P(\alpha)$ and $c_P^*(m_{\alpha})$.

\begin{theorem}\label{thm:dynamics}
For any $\alpha \in \C_k$ and $P\in \G_k$, the number $c_P^{*}(m_{\alpha})$ is divisible by $\frac{c_P(\alpha)}{\gcd(c_P(\alpha), k)}$ and divides $c_P(\alpha)$. In particular, for any $P\in \G_k$ and $\alpha\in \C_k$, the following hold:
\begin{enumerate}[(i)]
\item If $m_{\alpha}$ is fixed by $P$, then $c_P(\alpha)$ divides $k$, i.e., $P^{(k)}(\alpha)=\alpha$.
\item  The numbers $c_P(\alpha)$ and $c_P^*(m_{\alpha})$ satisfy the following inequality $$\frac{c_P(\alpha)}{k}\le c_P^{*}(m_{\alpha})\le  c_P(\alpha).$$ 
\item If $\gcd(c_{P}(\alpha), k)=1$, then $ c_P^{*}(m_{\alpha})=c_P(\alpha)$. In particular, if $\gcd(j, k)=1$ for any $j\in S_P$, then $S_P=S_P^*$.
\item The numbers $\mu_k^*(P)$ and $\mu_k(P)$ satisfy the following inequality $$\frac{\mu_k(P)}{k}\le \mu_k^*(P)\le \mu_k(P).$$
\end{enumerate}
\end{theorem}

\begin{proof}
For $\alpha\in \C_k$, set $i=c_P(\alpha)$ and $j=c_P^*(m_{\alpha})$. Since $P^{(i)}(\alpha)=\alpha$, from definition, $P^{(i)}\ast m_{\alpha}=m_{\alpha}$, hence $j$ divides $i$. Also, if $P^{(j)}\ast m_{\alpha}=m_{\alpha}$, then $P^{(j)}(\alpha)=\alpha^{q^s}$ for some $0\le s\le k-1$. Since $P$ has coefficients in $\F_q$ and $\alpha\in \F_{q^k}$, it follows that $P^{(jk)}(\alpha)=\alpha^{q^{sk}}=\alpha$, hence $jk$ is divisible by $i$. If we write $i_0=\gcd(i, k)$, we see that $ji_0$ is also divisible by $i$, and so $j$ is divisible by $\frac{i}{i_0}$. This shows that $c_P^{*}(m_{\alpha})$ is divisible by $\frac{c_P(\alpha)}{\gcd(c_P(\alpha), k)}$ and divides $c_P(\alpha)$. The proofs of (i), (ii), (iii) and $(iv)$ are straightforward.
\end{proof}

In particular, the previous theorem entails that the orbit length of a point $\alpha\in \C_k$ ``contracts'' in a factor at most $k$ to the orbit length of $m_{\alpha}\in \I_k$. As we further see, the bounds in item (ii) of Theorem~\ref{thm:dynamics} can be reached.

\begin{figure}[H]
  \centering {\includegraphics[width=1.0\linewidth]{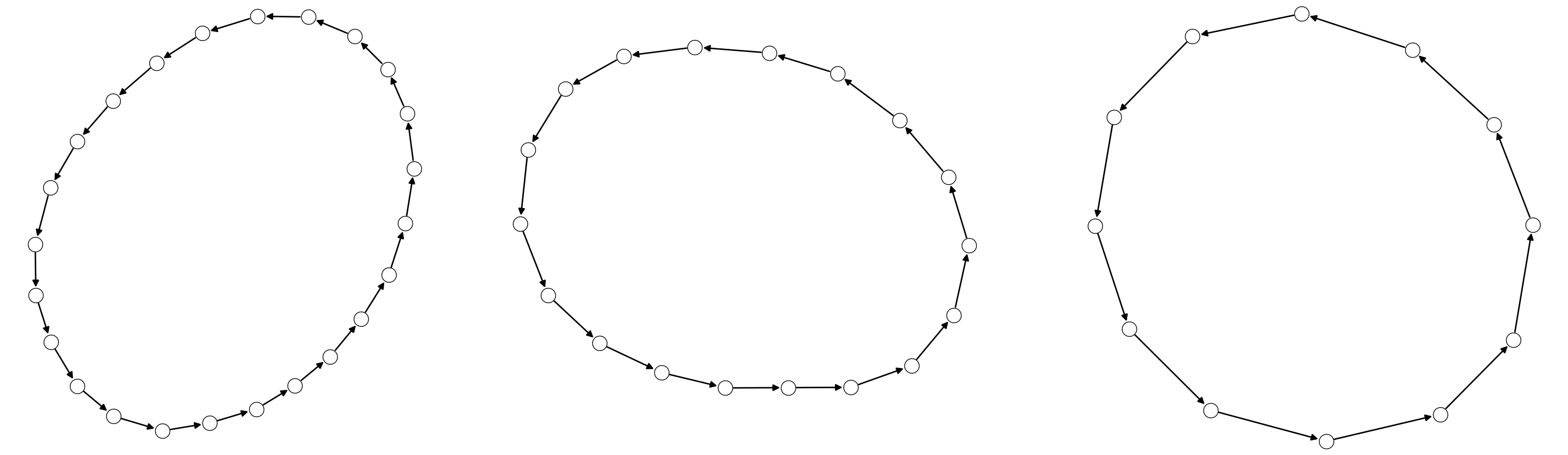}}
  {\caption{The functional graph $G(x^{13}+1, \C_6)$ for $q=2$.}
  \label{fig:dynamics}}
  \centering
\end{figure}

\begin{figure}[H]
  \centering
{\includegraphics[width=0.4\linewidth]{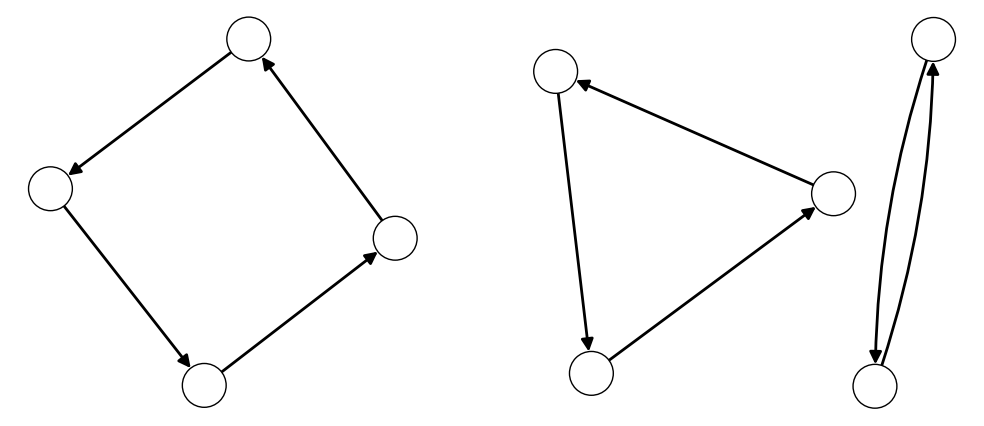}} 
 {\caption{The functional graph $G(x^{13}+1, \I_6)$ for $q=2$.}
  \label{fig:dynamics-2}}
  \centering
\end{figure}

\begin{example}
For $q=2$ and $k=6$, we observe that the polynomial $P=x^{13}+1$ permutes $\F_{64}=\F_{2^6}$. The functional graphs $G(P, \C_6)$ and $G(P, \I_6)$ are shown in Fig.~\ref{fig:dynamics} and~\ref{fig:dynamics-2}, respectively. In particular, we see that each cycle of $G(P, \C_6)$ contracts by the factor $6=k$ to a cycle of $G(P, \I_6)$.
\end{example}


Furthermore, we emphasize that the lower bound in item (ii) of Theorem~\ref{thm:dynamics} can always be attained. Fix $k\ge 1$ and consider $P=x^q$. Then $P\in \G_k$ and we can see that, for any $\alpha\in \C_k$, $c_P(\alpha)=k$. In addition, it follows from definition that $P\ast f=f$ for any $f\in \I_k$ and, in particular, $c_P^{*}(f)=1$ for any $f\in \I_k$. In other words, for any $\alpha\in \C_k$, the following holds
$$1=c_P^*(m_{\alpha})=\frac{c_P(\alpha)}{k}.$$

\subsection{The functional graph of certain permutations}
So far we have provided connections between the functional graphs of the maps induced by a permutation $P\in \G_k$ on the sets $\C_k$ and $\I_k$. The study of functional graphs of polynomial maps is usually made over the entire field: in \cite{A69}, monomial permutations are explored via the multiplicative order of finite fields and, in \cite{MV88}, the case of linearized permutations is considered, where the $\F_q$-order is employed. In general, there is no deterministic method to determine the functional graph of an arbitrary permutation (other than constructing the whole graph). The cases of monomial and linearized permutations are very special, since they are naturally connected to the algebraic structure of finite fields. 

Following the techniques employed in \cite{A69} and \cite{MV88}, we shall obtain the complete description of the dynamics of monomial and linearized permutations in the set $\C_k$. First, we fix some graph notation: for each positive integer $n$, $\cc(n)$ denotes the cyclic graph of length $n$. 

\begin{theorem}\label{thm:monomial-cycle}
Let $P=x^n$ with $\gcd(n, q^k-1)=1$ be a permutation of $\C_k$. The following holds:
\begin{equation*}
G(x^n, \C_k)=\bigoplus_{\ord_eq=k}\frac{\varphi(e)}{\ord_en}\times \cc(\ord_en).
\end{equation*}

In particular, $\mu_k^*(x^n)\ge \min\limits_{\ord_eq=k}\frac{\ord_e n}{k}$.
\end{theorem}

\begin{proof}
From Lemma~\ref{thm:mult}, $\C_k$ comprises the elements with multiplicative order $e$, where $e$ varies through the positive integers such that $k=\ord_eq$. In addition, for each positive integer $e$ with such property, $\C_k$ has $\varphi(e)$ elements with multiplicative order $e$. Therefore, for each positive integer $e$ such that $k=\ord_eq$, we just need to see how the elements with multiplicative order $e$ are distributed in the graph $G(x^n, \C_k)$. Let $\alpha\in \C_k$ be an element with multiplicative order $e$. From definition, $\alpha$ belongs to a cycle of length $j$ if and only if $j$ is the least positive integer such that $P^{(j)}(\alpha)=\alpha$, i.e., $\alpha^{n^j}=\alpha$. The latter is equivalent to $n^j\equiv 1\pmod e$. From definition, $j=\ord_en$. In particular, we have shown that each element of multiplicative order $e$ belongs to a cycle of length $j=\ord_en$. Since $\gcd(n, q^k-1)=1$, $\alpha$ and $\alpha^n$ have the same multiplicative order (see item (ii) of Lemma~\ref{thm:mult}) and, in particular, this shows that elements in a same cycle of $G(x^n, \C_k)$ have the same multiplicative order. Therefore, the divisor $e$ of $q^k-1$ contributes with $\frac{\varphi(e)}{\ord_en}$ copies of the cyclic graph of length $\ord_en$. The inequality $\mu_k^*(x^n)\ge \min_{\ord_eq=k}\frac{\ord_e n}{k}$ follows directly from item (iv) of Theorem~\ref{thm:dynamics}.
\end{proof}

Employing similar ideas (with suitable changes in the multiplicative-additive analogues), we obtain the linearized case.

\begin{theorem}\label{thm:linear-cycle}
Let $P=L_f$ be a linearized permutation of $\C_k$, where $L_f$ is the $q$-associate of $f\in \F_q[x]$ and $\gcd(f, x^k-1)=1$. The following holds:
\begin{equation}\label{eq:cycle-linear}
G(L_f, \C_k)=\bigoplus_{\ou(x, g)=k}\frac{\Phi_q(g)}{\ou(f, g)}\times \cc(\ou(f, g)),
\end{equation}
where $g\in \F_q[x]$ is monic. In particular, $\mu_k^*(L_f)\ge \min\limits_{\ou(x, g)=k}\frac{\ou(f, g)}{k}$.
\end{theorem}

\begin{proof}
From Theorem~\ref{thm:add}, $\C_k$ comprises the elements with $\F_q$-order $g$, where $g$ varies through the monic polynomials in $\F_q[x]$ such that $k=\ou(x, g)$. In addition, for each monic polynomial $g$ with such property, $\C_k$ has $\Phi_q(g)$ elements with $\F_q$-order $g$. Therefore, for each monic polynomial $g$ such that $k=\ou(x, g)$, we only need to show how the elements with $\F_q$-order $g$ are distributed in the graph $G(L_f, \C_k)$. Let $\alpha\in \C_k$ be an element with $\F_q$-order $g$. We observe that $\alpha$ belongs to a cycle of length $j$ if and only if $j$ is the least positive integer such that $P^{(j)}(\alpha)=\alpha$, i.e., $L_{f^j}(\alpha)=\alpha$. The latter is equivalent to $L_{f^j-1}(\alpha)=0$, i.e., $f^j\equiv 1\pmod g$. From definition, $j=\ou(f, g)$. In particular, we have shown that each element of $\F_q$-order $g$ belongs to a cycle of length $j=\ou(f, g)$. Since $\gcd(f, x^k-1)=1$, $\alpha$ and $L_f(\alpha)$ have the same multiplicative order (see item (ii) of Theorem~\ref{thm:add}) and, in particular, this shows that elements in a same cycle of $G(L_f, \C_k)$ have the same $\F_q$-order. Therefore, the monic divisor $g$ of $x^k-1$ contributes with $\frac{\Phi_q(g)}{\ou(f, g)}$ copies of the cyclic graph of length $\ou(f, g)$. Inequality $\mu_k^*(L_f)\ge \min\limits_{\ou(x, g)=k}\frac{\ou(f, g)}{k}$ follows directly from item (iv) of Theorem~\ref{thm:dynamics}.
\end{proof}

\subsubsection{M\"{o}bius dynamics}
When considering M\"{o}bius maps, the dynamics is quite trivial. For $A\in \GL_2(\F_q)$ with $A=\left(\begin{matrix}
a&b\\
c&d
\end{matrix}\right)$, let $\gamma_{A}:\F_{q^k}\to \F_{q^k}$ be the map given by $\gamma_{A}(\alpha)=\frac{a\alpha+b}{c\alpha+d}$ (with a suitable evaluation at the possible pole of $\gamma_A$). In Subsection~\ref{subsec:permut} we have seen that $\gamma_A$ induces a permutation on $\F_{q^k}$ and, since $A\in \GL_2(\F_q)$, such a permutation admits a polynomial representation $P_{A}\in \F_q[x]$. Therefore, from Proposition~\ref{local-permut}, $P_{A}$ permutes the set $\C_k$ (or, equivalently, $\gamma_A$ permutes $\C_k$). We observe that the possible pole of $\gamma_A$ is in $\F_q$. For simplicity, we assume $k\ge 2$ (the case $k=1$ can be easily studied in details). For $k\ge 2$, $\C_k$ has no pole of any M\"{o}bius map with coefficients in $\F_q$. In particular, we can iterate the function $\gamma_A$ on $\C_k$ without any consideration on the possible poles. From direct calculations, we have the equality of maps $\gamma_A^{(n)}=\gamma_{A^n}, n\in \mathbb N$. 
In particular, if $D=\ord([A])$ is the order of the class of $A$ in $\PGL_2(\F_q)$, $\gamma_A^{(D)}$ is the identity map. Therefore, $\gamma_A^{(D)}(\alpha)=\alpha$ for any $\alpha\in \C_k$ and so any point has period a divisor of $D$.

For $k\ge 3$, we observe that no point has period strictly smaller; if $\alpha\in \C_k$ with $k\ge 3$ and $\gamma_A^{(n)}(\alpha)=\alpha$, then $\gamma_{A^n}(\alpha)=\alpha$. If $[A^n]=[A]^n$ is not the identity $[I]$ of $\PGL_2(\F_q)$, the equality $\gamma_{A^n}(\alpha)=\alpha$ yields a nontrivial linear combination of $1, \alpha$ and $\alpha^2$ with coefficients in $\F_q$. This implies that $\alpha$ has minimal polynomial of degree at most $2$, a contradiction with $k\ge 3$. All in all, the previous observations easily imply the following result.

\begin{theorem}
Let $k\ge 3$ be a positive integer. For $A\in \GL_2(\F_q)$ with $A=\left(\begin{matrix}
a&b\\
c&d
\end{matrix}\right)$, the map $\gamma_{A}:\C_k\to \C_k$ given by $\gamma(\alpha)=\frac{a\alpha+b}{c\alpha+d}$ is well defined and is a permutation of $\C_k$. Additionally, if $n_k:= |\C_k|$ and $D$ denotes the order of $[A]$ in $\PGL_2(\F_q)$, the following holds
\begin{equation}
G(\gamma_A, \C_k)=\frac{n_k}{D}\times \cc(D).
\end{equation}
In other words, any element $\alpha\in \C_k$ is periodic with period $D$. 
\end{theorem}

In particular, we obtain the following corollaries.

\begin{corollary}\label{cor:mobius-divisible}
For $A\in \GL_2(\F_q)$ and $f\in \I_k$ with $k\ge 2$, let $\gamma_A\ast f$ be defined as in Subsection~\ref{subsec:Mobius-action}. If $D$ is the order of $[A]$ in $\PGL_2(\F_q)$ and $\gamma_A\ast f=f$, then $k=2$ or $k$ is divisible by $D$.
\end{corollary}

\begin{proof}
Let $\alpha\in \C_k$ be any root of $f$, hence $f=m_{\alpha}$ is the minimal polynomial of $\alpha$. From Theorem~\ref{thm:dynamics}, if $\gamma_A\ast f=f$ and $c_{P}(\alpha)$ denotes the least period of $\alpha$ by $\gamma_A$, then $c_P(\alpha)$ divides $k$. From the previous theorem, $c_P(\alpha)=D$ if $k\ge 3$. In particular, $k=2$ or $k$ is divisible by $D$.
\end{proof}
In particular, if $[A]$ has order $D$ in $\PGL_2(\F_q)$, the map $f\mapsto \gamma_A\ast f$ over $\I_k$ is free of fixed points whenever $k>2$ is not divisible by $D$.

\begin{corollary}
For $A\in \GL_2(\F_q)$ and $\alpha\in \C_k$ with $k\ge 2$, let $\gamma_A\ast m_{\alpha}$ be defined as in Subsection~\ref{subsec:Mobius-action}. If $D$ is the order of $[A]$ and $\gcd(D, k)=1$ then, for the permutation $P=\gamma_A$, we have that
$$c_{P}(\alpha)=c_P^*(m_{\alpha}).$$
\end{corollary}

In particular, we may attain the upper bound in item (ii) of Theorem~\ref{thm:dynamics}.

\section{Iterated construction of irreducible polynomials}
We have seen that if $f\in \F_q[x]$ is an irreducible polynomial of degree $k$ and $P\in \G_k$, then $$P\ast f=\gcd (f(P(x)), x^{q^k}-x),$$ is another irreducible polynomial of degree $k$. This identity suggests a recursive method for constructing irreducible polynomials of degree $k$ from a given $f\in \I_k$: we set $f_0=f$ and, for $i\ge 1$, $$f_i:= P\ast f_{i-1}=\gcd (f_{i-1}(P(x)), x^{q^k}-x).$$ From Definition~\ref{def:periods}, $j=c_{P}^*(f)$ is the least positive integer such that $f_j=f_0$: in particular, the sequence $\{f_i\}_{i\ge 0}$ has $c_P^*(f)$ distinct elements. We want to find a good permutation that generates a large number of irreducible polynomials from a single irreducible polynomial $f$. From the trivial bound $c_P^*(f)\ge \mu_k^*(P)$, it is sufficient to find permutations $P$ for which $\mu_k^*(P)$ is large. Since $\mu_k^*(P)\ge 1/k\cdot \mu_k(P)$, if $\mu_k(P)$ is large, then so is $\mu_k^*(P)$. We have the trivial bound $\mu_k^*(P)\le |\I_k|\approx \frac{q^k}{k}$. Theorem~\ref{thm:onto} shows that any permutation of $\I_k$ can be viewed as a map $f\mapsto P\ast f$ with $P\in \G_k$ and, in particular, there exist permutations $P\in \G_k$ for which $G(P, \I_k)$ comprises a \emph{full cycle} (i.e., all the elements of $\I_k$ lie in the same orbit): in this case, we have $\mu_k^*(P)=|\I_k|$. However, the construction of such permutations seems to be out of reach: in fact, even the construction of permutations of finite fields with a full cycle is not completely known \cite{CMT08}. Having this in mind, it would be interesting to obtain permutations $P$ for which $\mu_k^*(P)$ is reasonable large. 

As pointed out earlier, the description of the functional graph of general permutations of finite fields is still an open problem. However, for $P$ a monomial or a linearized polynomial, things are well understood and Theorems~\ref{thm:monomial-cycle} and \ref{thm:linear-cycle} provide lower bounds for the quantity $\mu_k^*(P)$: the lower bound for monomials depend on the positive integers $e$ for which $\ord_eq=k$ and, in the linearized case, the lower bound depends on the monic polynomials $g\in \F_q[x]$ for which $\ou(x, g)=k$. The numbers $e$ and the polynomials $g$ depend on the prime factorization of $q^k-1$ and the factorization of $x^k-1$ over $\F_q$, respectively.

In this section, we explore specific cases when $\mu_k^*(P)$ can be explicitly given and they are asymptotically the best possible in some cases.  We suppose that either $\frac{q^k-1}{q-1}$ is a prime number or $\frac{x^k-1}{x-1}$ is an irreducible polynomial. Though these conditions are very particular, they are reasonable natural: for instance, if $q=2$, then $\frac{q^k-1}{q-1}$ equals $2^k-1$ and the primes of this form are the \emph{Mersenne Primes}. In addition, if $k$ is a prime number, then $\frac{x^k-1}{x-1}=E_k(x)$ is just the $k$-th \emph{cyclotomic polynomial}: if $q$ is a primitive root modulo $k$, according to Theorem 2.47 of~\cite{LN}, $E_k(x)$ is an irreducible polynomial. We start with the monomial case.

\begin{proposition}
Suppose that $r=\frac{q^k-1}{q-1}$ is a prime number and let $P=x^n$ be a monomial permutation of $\F_{q^k}$ (i.e., $\gcd(n, q^k-1)=1$). For any polynomial $f\in \I_k$, let $\{f_i\}_{i\ge 0}$ be the sequence of polynomials given by $f_0=f$ and, for $i\ge 1$, $f_i:= P\ast f_{i-1}=\gcd (f_{i-1}(x^n), x^{q^k}-x)$. Then $\{f_i\}_{i\ge 0}$ yields at least $\frac{\ord_rn}{k}$ irreducible polynomials of degree $k$. In particular, if $n$ is a primitive root modulo $r$, then $\frac{\ord_rn}{k}=\frac{q^k-q}{(q-1)k}$.
\end{proposition}

\begin{remark}
We observe that if $r=\frac{q^k-1}{q-1}$ is prime, then $k$ is a prime number. In this case, $|\I_k|=\frac{q^k-q}{k}$ and so, if $n$ is a primitive root modulo $r$ and $f$ is any element of $\I_k$, the sequence $\{f_i\}_{i\ge 0}$ produces at least $\frac{1}{q-1}$ 
of the elements in $\I_k$. 
\end{remark}

We now proceed to the linearized case.

\begin{proposition}\label{prop:cycle-cyclotomic}
Suppose that $E_k(x)=\frac{x^k-1}{x-1}$ is an irreducible polynomial and let $P=L_g\in \F_q[x]$ be a linearized permutation of $\F_{q^k}$ (i.e., $\gcd(g, x^k-1)=1$). For any polynomial $f\in \I_k$, let $\{f_i\}_{i\ge 0}$ be the sequence of polynomials given by $f_0=f$ and, for $i\ge 1$, $f_i:= P\ast f_{i-1}=\gcd (f_{i-1}(L_g(x)), x^{q^k}-x)$. Then $\{f_i\}_{i\ge 0}$ yields at least $\frac{\ou(g, E_k)}{k}$ irreducible polynomials of degree $k$. In particular, if $\ou(g, E_k)=\Phi_q(E_k)=q^{k-1}-1$, then $\frac{\ou(g, E_k)}{k}=\frac{q^{k-1}-1}{k}$.
\end{proposition}

\begin{remark}
We observe that if $E_k(x)=\frac{x^k-1}{x-1}$ is an irreducible polynomial, then $k$ is a prime number. In this case, $|\I_k|=\frac{q^k-q}{k}$ and so, if $\ou(f, E_k)=\Phi_q(E_k)=q^{k-1}-1$ and $f$ is any element of $\I_k$, then the sequence $\{f_i\}_{i\ge 0}$ produces at least $\frac{1}{q}$ of the elements in $\I_k$. 
\end{remark}

If $E_k(x)$ is irreducible, then the quotient $K=\frac{\F_q[x]}{\langle E_k(x)\rangle}$ is a field which is isomorphic to $\F_{q^{k-1}}$. In this case, if $\theta$ denotes the class of $x$ in the quotient $\frac{\F_q[x]}{\langle E_k(x)\rangle}$, it is direct to verify that $\ou(g, E_k)=\ord(g(\theta))$. Therefore, $\ou(g, E_k)=\Phi_q(E_k)=q^{k-1}-1=|\F_{q^{k-1}}^*|$ if and only if $g(\theta)$ is a \emph{primitive} element of $\F_{q^{k-1}}$. 

Primitive elements play important roles in cryptography and coding theory and have been extensively studied in the past few decades. However, the efficient construction of primitive element in finite fields $\F_{q^n}$ is still an open problem. Mainly, this is due to the fact that a general method for constructing such elements requires the prime factorization of $q^n-1$. Nevertheless, many authors have been treated the problem to find elements with a reasonable \emph{high multiplicative order} in finite fields as a relaxation of primitive elements. High order elements have been used in many practical ways, including cryptography, pseudo random number generator and the construction of gauss periods \cite{ASV}. Most notably, high order elements were employed in the celebrated AKS primality test~\cite{AKS}. Due to these applications, the construction of high order elements have been considered by many authors in the past years. In \cite{G99}, this problem is treated in general extensions of finite fields and some particular extensions are further more explored, including Artin-Schreier extensions~\cite{BR, V04} and cyclotomic extensions~\cite{P14}. The latter covers exactly the kind of extensions that we are considering here: from Theorem 2 of \cite{P14}, we obtain the following result.

\begin{proposition}\label{prop:popovych}
Let $q$ be a power of a prime $p$. If $a$ is any nonzero element of $\F_q$ and $q$ is primitive modulo the prime $k$, then the class of $x+a$ in $K=\frac{\F_q[x]}{\langle E_k(x)\rangle}$ has multiplicative order at least $\tau(p, k)$, where
$\tau(2, k)=2^{\sqrt{2(k-2)}-2}$, $\tau(3, k)=3^{\sqrt{3(k-2)}-2}$ and $\tau(p, k)=5^{\sqrt{(k-2)/2}-2}$ for $p\ge 5$.
\end{proposition}

We observe that any nonzero element of $\F_{q^{k-1}}$ is written as $h(\theta)$, for some polynomial $h$ of degree at most $k-1$ such that $h$ is not divisible by $E_k$. In particular, if we know that $h(\theta)\ne 0$ has multiplicative order $e=\ord(h(\theta))$ in $\F_{q^{k-1}}$, we have a method to produce at least $\frac{e}{k}$ irreducible polynomials of degree $k$ from a single $f\in \I_k$: since $h(\theta)\ne 0$, $h$ is a polynomial of degree at most $k-1$ such that $\gcd(h, E_k)=1$. In particular, there exists a polynomial $H$ of degree at most $k$ such that $\gcd(H(x), x-1)=1$ and $H\equiv h\pmod {E_k}$: in this case, $\gcd(H(x), x^k-1)=1$ and so $L_H$, the $q$-associate of $H$, is a linearized permutation of $\F_{q^k}$. If $f$ is any irreducible polynomial of degree $k$, it follows from Proposition~\ref{prop:cycle-cyclotomic} that the sequence $\{f_i\}_{i\ge 0}$ given by $f_0=f$ and $f_i=\gcd(f_{i-1}(L_H(x)), x^{q^k}-x)$ yields at least $\frac{\ou(g, E_k)}{k}=\frac{\ord(h(\theta))}{k}=\frac{e}{k}$ distinct irreducible polynomials of degree $k$. 

Using this approach, Proposition~\ref{prop:popovych} suggests to consider $h(x)=x+a$, where $a$ is a nonzero element of $\F_q$: in this case, the same proposition provides the bound $e\ge \tau(p, k)$. For $h(x)=x+a$, $H$ is a polynomial of degree at most $k$ such that $H\equiv x+a\pmod {E_k}$ and $\gcd(H(x), x-1)=1$. If $q\ne 2$, there exists $a\in \F_q\setminus \{0, 1\}$ and so $H(x)=x-a$ satisfies the required properties. For $q=2$, $a=1$ is the only nonzero element of $\F_2$ and so we consider $H(x)=E_k(x)+x+a=E_k(x)+x+1$: we observe that, since $q=2$ is a primitive root modulo $k$, $k$ is not $2$ and so $H(1)=E_k(1)+1+1=k+2=k$ which is not zero in $\F_q$, since $k\ne 2$ is a prime. Therefore, for $q=2$, we use $H(x)=E_k(x)+x+1=x^{k-1}+\cdots+x^2$. In summary,  we obtain the following result.

\begin{theorem}
Let $q$ be a power of a prime $p$ and let $k$ be a prime number such that $q$ is a primitive root modulo $k$. In addition, let $\tau(p, k)$ be the function given in Proposition~\ref{prop:popovych}. Then, for any $f\in \I_k$, the sequence $\{f_i\}_{i\ge 0}$ given by $f_0=f$ and $f_i(x)=\gcd (f_{i-1}(L_H(x)), x^{q^k}-x)$ if $i\ge 1$ yields at least $\frac{\tau(p, k)}{k}$ distinct irreducible polynomials of degree $k$ in each of the following cases:
\begin{enumerate}[(i)]
\item $q=2$ and $L_H(x)=x^{q^{k-1}}+\cdots+x^{q^2}=x^{2^{k-1}}+\cdots+x^{4}$,
\item $q\ne 2$ and $L_H(x)=x^{q}-ax$ with $a\ne 0, 1$ and $a\in \F_q$.
\end{enumerate}
\end{theorem}

The following example provides a numerical instance of the previous theorem.
\begin{example}
We observe that $k=53$ is prime and $3$ is a primitive root modulo $53$. In particular, for any irreducible polynomial $f\in \F_3[x]$ of degree $53$, the sequence $\{f_i\}_{i\ge 0}$ given by $f_0=f$ and $f_i(x)=\gcd(f_{i-1}(x^3+x), x^{3^{53}}-x)$ if $i\ge 1$ yields at least
$$\left\lceil\frac{\tau(3, 53)}{53}\right\rceil=1,672,$$
distinct irreducible polynomials of degree $53$ over $\F_3$. For instance, one may consider the initial input $f_0(x)=x^{53} - x^4 -x^3 -x^2 + 1\in \F_3[x]$. We emphasize that the acutal order of the class of $x+1$ in the quotient $K=\F_3[x]/(\Phi_{53}(x))=\F_{3^{53}}$ is $134718888901384\approx 1.3\cdot 10^{14}$ and so 
$$\left\lceil\frac{134718888901384}{53}\right\rceil \approx 2.5\times 10^{12}.$$
\end{example}

\end{document}